\documentclass{amsart}

\usepackage{amssymb}
\usepackage{stmaryrd}
\usepackage{mathrsfs}
\usepackage{array,subfigure}
\usepackage{pictex}

\newtheorem{thm}{Theorem}[section]
\newtheorem{prop}[thm]{Proposition}
\newtheorem{cor}[thm]{Corollary}

\newtheorem{defn}[thm]{Definition}
\newenvironment{xpl}{\medskip \noindent {\bf  Example.}}{\hfill$\diamondsuit$\mbox{}\bigskip}
\newcounter{num}
\newenvironment{thmlist}{\begin{list}{(\roman{num})}{\usecounter{num}\setlength{\leftmargin}{25pt}
\setlength{\itemindent}{0pt}\setlength{\labelwidth}{20pt}\setlength{\labelsep}{5pt}\setlength{\itemsep}{0in}
\setlength{\listparindent}{0pt}}}{\end{list}}
\numberwithin{equation}{section}

\renewcommand{\appendix}{\renewcommand{\section}{
    \secdef\Appendix\sAppendix}
\setcounter{section}{0}
\renewcommand{\thesection}{\Alph{section}}}

\newcommand{\Appendix}[2][default]{
\refstepcounter{section}
\addcontentsline{toc}{section}{\protect\numberline{\appendixname~\thesection}{\hspace{50pt} #1}}
{\flushleft\Large\bfseries\appendixname\ \thesection\hspace{1em}#2\par}}

\newcommand{\sAppendix}[1]{
{\flushleft\large\bfseries\appendixname\par\flushleft#1\par}
\vspace{\baselineskip}}

\newcommand{\C}{\mathbb{C}}
\newcommand{\Ha}{\mathbb{H}}
\newcommand{\R}{\mathbb{R}}
\newcommand{\Z}{\mathbb{Z}}
\newcommand{\Q}{\mathbb{Q}}
\newcommand{\cps}{\mathbb{C}P}
\newcommand{\qps}{\mathbb{H}P}
\newcommand{\ol}[1]{\bar{#1}}
\newcommand{\sm}[1]{\scriptscriptstyle{#1}}

\newcommand{\Isom}{\operatorname{Isom}}

\newcommand{\Aut}{\operatorname{Aut}}
\newcommand{\Hom}{\operatorname{Hom}}
\newcommand{\im}{\operatorname{Im}}
\newcommand{\lcm}{\operatorname{lcm}}

\newcommand{\ric}{\operatorname{Ricci}}
\newcommand{\Ric}{\operatorname{Ric}}
\newcommand{\ord}{\operatorname{Ord}}
\newcommand{\ind}{\operatorname{Ind}}
\newcommand{\Vol}{\operatorname{Vol}}

\newcommand{\Spec}{\operatorname{Spec}}
\newcommand{\SF}{\operatorname{SF}}

\newcommand{\picorb}{\operatorname{Pic^{\text{orb}}}}

\newcommand{\Div}{\operatorname{div}}
\newcommand{\divorb}{\operatorname{Div^{\text{orb}}}}
\newcommand{\weil}{\operatorname{Weil}}
\newcommand{\orb}{\operatorname{Orb}}
\newcommand{\Sing}{\operatorname{Sing}}

\begin{document}
\pagenumbering{roman}

\title{Some Examples of toric Sasaki-Einstein Manifolds}
\author{Craig van Coevering}
\address{Department of Mathematics, Massachusetts Institute of Technology, 77 Massachusetts Avenue, Cambridge, MA 02139-4307}
\email{craig@math.mit.edu}
\date{January 18, 2007}

\begin{abstract}
A series of examples of toric Sasaki-Einstein 5-manifolds is constructed which first appeared in
the author's Ph.D. thesis~\cite{vC2}.  These are submanifolds of the toric 3-Sasakian 7-manifolds of C. Boyer and K. Galicki.
And there is a unique toric quasi-regular Sasaki-Einstein 5-manifold associated to every simply
connected toric 3-Sasakian 7-manifold.  Using 3-Sasakian reduction as in~\cite{BGMR,BGM2} an infinite series of examples is
constructed of each odd second Betti number.  They are all diffeomorphic to $\#k M_\infty$, where
$M_\infty \cong S^2\times S^3$, for $k$ odd.
We then make use of the same framework to construct positive
Ricci curvature toric Sasakian metrics on the manifolds $X_\infty \# kM_\infty$ appearing in the classification
of simply connected smooth 5-manifolds due to Smale and Barden.  These manifolds are not spin, thus do not
admit Sasaki-Einstein metrics.  They are already known to admit toric Sasakian metrics (cf.~\cite{BGO})
which are not of positive Ricci curvature.
We then make use of the join construction of C. Boyer and K. Galicki first appearing in~\cite{BG3}, see also~\cite{BGO},
to construct infinitely many toric Sasaki-Einstein manifolds with arbitrarily high second Betti number
of every dimension $2m+1\geq 5$.  This is in stark contrast to the analogous case of Fano manifolds in even dimensions.
\end{abstract}

\maketitle

\pagestyle{plain}



\pagenumbering{arabic}

\section{Introduction}

A new series of quasi-regular Sasaki-Einstein 5-manifolds is constructed.  These examples first appeared
in the author's Ph.D. thesis~\cite{vC2}.  They are toric, and arise as submanifolds of toric 3-Sasakian
7-manifolds.  Applying 3-Sasakian reduction to torus actions on spheres C. Boyer, K. Galicki,
\textit{et al}~\cite{BGMR} produced infinitely many toric 3-Sasakian 7-manifolds.  More precisely,
one has a 3-Sasakian 7-manifold $\mathcal{S}_\Omega$ for each integral weight matrix $\Omega$ satisfying some
conditions to ensure smoothness.  This produces infinitely many examples of each $b_2(\mathcal{S}_\Omega)\geq 1$.
A result of D. Calderbank and M. Singer~\cite{CaSin} shows that, up to finite coverings, this produces
all examples of toric 3-Sasakian 7-manifolds.  Associated to each $\mathcal{S}_\Omega$ is its twistor space
$\mathcal{Z}$, a complex contact Fano 3-fold with orbifold singularities.  The action of $T^2$ complexifies
to $T^2_{\C}=\C^*\times\C^*$ acting on $\mathcal{Z}$.  Furthermore, if $\mathbf{L}$ is the line bundle of
the complex contact structure, the action defines a pencil
\[E=\mathbb{P}(\mathfrak{t}_\C)\subseteq|\mathbf{L}|,\]
where $\mathfrak{t}_\C$ is the Lie algebra of $T^2_{\C}$.  We determine the structure of the divisors in $E$.
The generic $X\in E$ is a toric variety with orbifold structure whose orbifold anti-canonical bundle
$\mathbf{K}^{-1}_X$ is positive.  The total space $M$ of the associated $S^1$ orbifold bundle to $\mathbf{K}_X$
has a natural Sasaki-Einstein structure.  Associated to any toric 3-Sasakian 7-manifold $\mathcal{S}$ with
$\pi_1(\mathcal{S})=e$ we have the following diagram.
\begin{equation}\label{int:diag-cor}
\beginpicture
\setcoordinatesystem units <1pt, 1pt> point at 0 30
\put {$M$} at -15 60
\put {$\mathcal{S}$} at 15 60
\put {$X$} at -15 30
\put {$\mathcal{Z}$} at 15 30
\put {$\mathcal{M}$} at 15 0
\arrow <2pt> [.3, 1] from -6 60 to 6 60
\arrow <2pt> [.3, 1] from -6 30 to 6 30
\arrow <2pt> [.3, 1] from -15 51 to -15 39
\arrow <2pt> [.3, 1] from 15 51 to 15 39
\arrow <2pt> [.3, 1] from 15 21 to 15 9
\endpicture
\end{equation}
The horizontal maps are inclusions and the vertical are orbifold fibrations.  And $\mathcal{M}$ is the
4-dimensional anti-self-dual Einstein orbifold over which $\mathcal{S}$ is an $Sp(1)$ or $SO(3)$ orbifold
bundle.

It follows from the Smale/Barden classification of smooth 5-manifolds (cf.~\cite{Sm} and~\cite{Bar}) that
$M$ is diffeomorphic to $\#k(S^2\times S^3)$ where $b_2(M)=k$.
We have the following:
\begin{thm}\label{int:main-thm}
Associated to every simply connected toric 3-Sasakian 7-manifold $\mathcal{S}$ is a toric quasi-regular
Sasaki-Einstein 5-manifold $M$.  If $b_2(\mathcal{S})=k$, then $b_2(M)=2k+1$ and
\[M\underset{\text{diff}}{\cong}\# m(S^2\times S^3), \text{ where  }m=b_2(M).\]
This gives an invertible correspondence.  That is, given either $X$ or $M$ in diagram~\ref{int:diag-cor}
one can recover the other spaces with their respective geometries.
\end{thm}
This gives an infinite series of quasi-regular Sasaki-Einstein structures on $\# m(S^2\times S^3)$
for every odd $m\geq 3$.

In section~\ref{sec:Sasak} we review the basics of Sasakian geometry.  In section~\ref{sec:toric}
we cover the toric geometry used in the proof of~\ref{int:main-thm}.  The reader can
find a complete proof of the solution to the Einstein equations giving theorem~\ref{int:main-thm} in~\cite{vC2}.
But the existence problem of Sasaki-Einstein structures on toric Sasakian manifolds is completely solved
in~\cite{WaZhu} and~\cite{FOW}.  So this result is summarized in section~\ref{sec:Einst-eq}.
The basics of 3-Sasakian manifolds and 3-Sasakian reduction are covered in section~\ref{sec:3Sasak}.
This section also contains some results on anti-self-dual Einstein orbifolds such as
the classification result of D. Calderbank and M. Singer used in the correspondence in theorem~\ref{int:main-thm}.
Section~\ref{sec:Sasak-submfd} contains the proof of theorem~\ref{int:main-thm} and completes diagram
~\ref{int:diag-cor}.  The only other possible topological types for simply connected toric Sasakian
5-manifolds are the non-spin manifolds $X_\infty \# kM_\infty$.  In section~\ref{sec:pos-Sasak} we use the framework
already constructed to construct positive Ricci curvature Sasakian metrics on these manifolds.
They are already shown to admit Sasakian structures in~\cite{BGO}.  In section~\ref{sec:high-dim}
we use the join construction of C. Boyer and K. Galicki~\cite{BG3,BGO} to construct higher dimensional examples of
toric quasi-regular Sasaki-Einstein and positive Ricci curvature manifolds.  In particular, in every
possible dimension $n\geq 5$ there exist infinitely many toric quasi-regular Sasaki-Einstein manifolds
with arbitrarily high second Betti number.

This article concentrates on the quasi-regular case.  This is not for lack of interest in irregular
Sasaki-Einstein manifolds, but that case is covered well elsewhere, such as in~\cite{MSY1,MSY2,FOW}.

This article makes copious use of orbifolds, V-bundles on them, and orbifold invariants such as the orbifold fundamental group
$\pi_1^{orb}$ and orbifold cohomology $H_{orb}^*$.  We will use the terminology 'V-bundle' to denote an orbifold bundle.
The characteristic classes of V-bundles will be elements of orbifold cohomology.
The reader unfamiliar with these concepts can consult~\cite{HS1} or
the appendices of~\cite{BG2} or~\cite{vC2}.

\section{Sasakian manifolds}\label{sec:Sasak}

We summarize the basics of Sasakian geometry in this section.  See the survey article of C. Boyer and K. Galicki~\cite{BG2}
for more details.
\begin{defn}\label{defn:Sasak}
Let $(M,g)$ be a Riemannian manifold of dimension $n=2m+1$, and $\nabla$ the Levi-Civita connection.
Then $(M,g)$ is Sasakian if either of the following
equivalent conditions hold:
\begin{thmlist}
\item There exists a unit length Killing vector field $\xi$ on $M$ so that the $(1,1)$ tensor
$\Phi(X)=\nabla_X\xi$ satisfies the condition
\[(\nabla_X\Phi)(Y)=g(\xi,Y)X-g(X,Y)\xi \]
for vector fields $X$ and $Y$ on $M$.
\item The metric cone $C(M)=\R_+\times M, \ol{g}=dr^2 +r^2g$ is K\"{a}hler.
\end{thmlist}
\end{defn}

Define $\eta$ to be the one form dual to $\xi$, i.e. $\eta(X)=g(X,\xi)$.
We say that $\{g,\Phi,\xi,\eta\}$ defines a Sasakian structure on $M$.
Note that $D=\ker\eta$ defines a contact structure on $M$, and $\Phi$ defines a
$CR$-structure on $D$.  Also, the integral curves of $\xi$ are geodesics.
The one form $\eta$ extends to a one form on $C(M)$ as $\eta(X)=\frac{1}{r^2}\ol{g}(\xi,X)$.
In (\textit{ii}) of the definition $M$ is identified with the subset $r=1$ of $C(M)$, and $\xi=Jr\partial_r$.
And the complex structure arises as
\begin{equation}
Jr\partial_r =\xi \quad JY =\Phi(Y)-\eta(Y)r\partial_r \quad\text{for } Y\in TM.
\end{equation}
The K\"{a}hler form of $(C(M),\ol{g})$ is given by
\begin{equation}
\frac{1}{2}d(r^2\eta)=\frac{1}{2}dd^c r^2,
\end{equation}
where $d^c =\frac{i}{2}(\ol{\partial} -\partial)$.

Besides the K\"{a}hler structure on $C(M)$ there is \emph{transverse K\"{a}hler structure} on $M$.
The Killing vector field $\xi$ generates the Reeb foliation $\mathscr{F}_\xi$ on $M$.
The vector field $\xi-iJ\xi=\xi+ir\partial_r$ on $C(M)$ is holomorphic and generates a local
$\C^*$-action extending that of $\xi$ on $M$.  The local orbits of this action define a transverse
holomorphic structure on $\mathscr{F}_\xi$.  One can choose an open covering $\{U_\alpha\}_{\alpha\in A}$
of $M$ such that we have the projection onto the local leaf space
$\pi_\alpha:U_\alpha \rightarrow V_\alpha$.  Then when $U_\alpha\cap U_\beta\neq\emptyset$, the transition
\[\pi_\beta\circ\pi_\alpha^{-1}:\pi_\alpha(U_\alpha\cap U_\beta)\rightarrow\pi_\beta(U_\alpha\cap U_\beta)\]
is holomorphic.  There is an isomorphism
\[d\pi_\alpha :D_p \rightarrow T_{\pi_\alpha (p)}V_\alpha,\]
for each $p\in U_\alpha$, which allows one to define a metric $g^T$ with K\"{a}hler form $\omega=\frac{1}{2}d\eta$
as restrictions of $g$ and $\frac{1}{2}d\eta=g(\Phi\cdot,\cdot)$ to $D_p$.  These are easily seen to be invariant
under coordinate changes.
Straight forward calculation gives the following:
\begin{prop}
Let $(M,g)$ be a Sasakian manifold, and $\pi_\alpha:U_\alpha \rightarrow V_\alpha$ as above.  If
$Y,Z\in\Gamma(D)$ are $\pi_\alpha$-related to $\tilde{Y},\tilde{Z}\in\Gamma(TV_\alpha)$ then
\begin{gather}
\Ric^T(\tilde{Y},\tilde{Z})=\Ric(Y,Z) +2g(X,Y),\label{eq:Sasak-Ric}\\
s^T = s+2m,
\end{gather}
where $\Ric, s$, resp. $\Ric^T, s^T$, are the Ricci and scalar curvatures of $g$, resp. $g^T$.
\end{prop}

\begin{defn}
A \emph{Sasaki-Einstein} manifold is a Sasakian manifold $(M,g)$ with $\Ric= 2mg$.
\end{defn}
Note that one always has $\Ric(\xi,\xi)=2m$ which fixes the Einstein constant at $2m$.  Simple calculation
shows that $(M,g)$ is Einstein if, only if, $(C(M),\ol{g})$ is Ricci flat.

If $\xi$ on $M$ induces a free $S^1$-action then $\{g,\Phi,\xi,\eta\}$ is a \emph{regular} Sasakian structure.
Another possibility is that all the orbits close but the action is not free, then the structure is \emph{quasi-regular}.
The third possibility is that not all the orbits close, in which case the generic orbit does not close.  In this case the
Sasakian structure is \emph{irregular}.
In the regular, resp. quasi-regular, cases the leaf space, along with its transverse K\"{a}hler structure,
is a manifold, resp. orbifold, $X$.  And $M$ is the total space of a principal $S^1$-bundle, resp. $S^1$ V-bundle,
$\pi:M\rightarrow X$.  A V-bundle is the orbifold analogue of a fiber bundle.  See the appendices of~\cite{BG2} or
~\cite{vC2} for details.

For a quasi-regular Sasakian manifold $(M,g,\Phi,\xi,\eta)$, the leaf space $X$ of $\mathscr{F}_\xi$ is a normal
projective, algebraic variety with an orbifold structure and a K\"{a}hler form $\omega$ with
$[\omega]\in H_{orb}^2(X,\Z)$.  We will make use of the following well known converse~\cite{BG3}.

\begin{prop}\label{KE-Sasak}
Let $(X,\omega)$ be a K\"{a}hler orbifold, with $[\omega]\in H^2_{orb}(X,\Z)$.  There is a holomorphic line
V-bundle $\mathbf{L}$ with $c_1(\mathbf{L})=-[\omega]\in H^2_{orb}(X,\Z)$ and an $S^1$-principal subbundle
$M\subset\mathbf{L}$ such that $M$ has a family of Sasakian structures $\{g_a,\Phi,\xi_a,\eta_a\}, a\in\R_+$,
\[g_a =a^2\eta\otimes\eta +a\pi^*h, \text{  where $h$ is a K\"{a}hler metric on }X. \]
Furthermore, if $\mathbf{L}^\times$ is $\mathbf{L}$ minus the zero section, then $\mathbf{L}^\times$ is biholomorphic to $C(M)$.

If $(X,\omega)$ is K\"{a}hler-Einstein with positive scalar curvature, and $q[\omega]\in c_1^{orb}(X)$, for $q\in\Z_+$,
then exactly one of the above Sasakian structures is Sasaki-Einstein, and one may take $\mathbf{L}=K_X^{\frac{1}{q}}$.
\end{prop}
Note that if $\pi_1^{orb}(X)=e$ and $[\omega]$ is indivisible in $H^2_{orb}(X,\Z)$, then $\pi_1^{orb}(M)=e$.
Of course we are interested in the case where $M$ is smooth.  This happens when the action of the local uniformizing groups
of $X$ inject into the group of the fibers of the V-bundle $M$.

We will mainly be interested in \emph{toric} Sasaki-Einstein manifolds.
\begin{defn}
A \emph{toric Sasakian manifold} is a Sasakian manifold $M$ of dimension $2m+1$ whose Sasakian structure
$\{g,\Phi,\xi,\eta\}$ is preserved by an effective action of an $(m+1)$-dimensional torus $T$ such that
$\xi$ is an element of the Lie algebra $\mathfrak{t}$ of $T$.
\end{defn}
Let $T_\C \cong(\C^*)^{m+1}$ be the complexification of $T$, then $T_\C$ acts on $C(M)$ by holomorphic automorphisms.
One sees that this definition is equivalent to $C(M)$ being a toric K\"{a}hler manifold.

\section{Toric geometry}\label{sec:toric}

We give some basic definitions in the theory of toric varieties that we will need.
See~\cite{Ful,O1,O2} for more details.  In addition we will consider the notion of a compatible
orbifold structure on a toric variety and holomorphic V-bundles.  We are interested in
K\"{a}hler toric orbifolds, and will give a description of the K\"{a}hler structure due to
V. Guillemin~\cite{Gu1}.

\subsection{Toric varieties}

Let $N\cong\Z^r$ be the free $\Z$-module of rank r and $M=\Hom_\Z(N,\Z)$ its dual.
We denote $N_\Q =N\otimes\Q$ and $M_\Q=M\otimes\Q$ with the natural pairing
\[\langle\ \, ,\ \rangle : M_\Q \times N_\Q \rightarrow\Q.\]
Similarly we denote $N_\R =N\otimes\R$ and $M_\R =M\otimes\R$.

Let $T_\C:=N\otimes_{\Z}\C^*\cong\C^*\times\cdots\times\C^*$ be the algebraic torus.
Each $m\in M$ defines a character $\chi^m :T_\C \rightarrow\C^*$ and each $n\in N$ defines a
one-parameter subgroup $\lambda_n :\C^*\rightarrow T_\C$.  In fact, this gives an isomorphism between
$M$ (resp. $N$) and the multiplicative group $\Hom_{\text{alg.}}(T_\C,\C^*)$
(resp. $\Hom_{\text{alg.}}(\C^*,T_\C)$).

\begin{defn}
A subset $\sigma$ of $N_\R$ is a \emph{strongly convex rational polyhedral cone} if there are
$n_1,\ldots, n_r$ so that
\[ \sigma=\R_{\geq 0}n_1 +\cdots +\R_{\geq 0}n_r,\]
and one has $\sigma\cap -\sigma=\{o\}$, where $o\in N$ is the origin.
\end{defn}
The dimension $\dim\sigma$ is the dimension of the $\R$-subspace $\sigma +(-\sigma)$ of $N_\R$.
The dual cone to $\sigma$ is
\[\sigma^{\vee}=\{x\in M_\R : \langle x,y\rangle\geq 0\text{ for all }y\in\sigma\},\]
which is also a convex rational polyhedral cone.
A subset $\tau$ of $\sigma$ is a face, $\tau <\sigma$, if
\[ \tau=\sigma\cap m^{\perp}=\{y\in\sigma : \langle m,y\rangle =0\}\text{ for }m\in\sigma^{\vee}.\]
And $\tau$ is a strongly convex rational polyhedral cone.
\begin{defn}
A \emph{fan} in $N$ is a collection $\Delta$ of strongly convex rational polyhedral cones
such that:
\begin{thmlist}
\item  For $\sigma\in\Delta$ every face of $\sigma$ is contained in $\Delta$.
\item  For any $\sigma,\tau\in\Delta$, the intersection $\sigma\cap\tau$ is a face of both $\sigma$
and $\tau$.
\end{thmlist}
\end{defn}
We will consider \emph{complete} fans for which the support $\bigcup_{\sigma\in\Delta}\sigma$ is $N_\R$.
We will denote
\[ \Delta(i):= \{\sigma\in\Delta : \dim\sigma =i\},\quad 0\leq i\leq n. \]
\begin{defn}
A fan in $N$ is \emph{nonsingular} if each $\sigma\in\Delta(r)$ is generate by $r$ elements of $N$
which can be completed to a $\Z$-basis of $N$.  A fan in $N$ is \emph{simplicial} if each $\sigma\in\Delta(r)$ is generated
by $r$ elements of $N$ which can be completed to a $\Q$-basis of $N_\Q$.
\end{defn}
If $\sigma$ is a strongly convex rational polyhedral cone, $S_\sigma = \sigma^{\vee}\cap M$ is a finitely generated semigroup.
We denote by $\C[S_\sigma]$ the semigroup algebra.  We will denote the generators of $\C[S_\sigma]$ by
$x^m$ for $m\in S_\sigma$.
Then $U_\sigma := \Spec\C[S_\sigma]$ is a normal affine variety
on which $T_\C$ acts algebraically with a (Zariski) open orbit isomorphic to $T_\C$.  If $\sigma$ is nonsingular, then
$U_\sigma \cong\C^n$.
\begin{thm}[\cite{Ful,O1,O2}]
For a fan $\Delta$ in $N$ the affine varieties $U_\sigma$ for $\sigma\in\Delta$ glue together to
form an irreducible normal algebraic variety
\[ X_\Delta =\bigcup_{\sigma\in\Delta}U_\sigma. \]
Furthermore, $X_\Delta$ is non-singular if, and only if, $\Delta$ is nonsingular.  And $X_\Delta$ is
compact if, and only if, $\Delta$ is complete.
\end{thm}
\begin{prop}\label{prop:toric-prop}
The variety $X_\Delta$ has an algebraic action of $T_\C$ with the following properties.
\begin{thmlist}
\item  To each $\sigma\in\Delta(i), 0\leq i \leq n,$ there corresponds a unique $(n-i)$-dimensional
$T_\C$-orbit $\orb(\sigma)$ so that $X_\Delta$ decomposes into the disjoint union
\[ X_\Delta =\bigcup_{\sigma\in\Delta}\orb(\sigma), \]
where $\orb(o)$ is the unique $n$-dimensional orbit and is isomorphic to $T_\C$.
\item  The closure $V(\sigma)$ of $\orb(\sigma)$ in $X_\Delta$ is an irreducible $(n-i)$-dimensional
$T_\C$-stable subvariety and
\[ V(\sigma)= \bigcup_{\tau\geq\sigma}\orb(\tau). \]
\end{thmlist}
\end{prop}

We will consider toric varieties with an orbifold structure.
\begin{defn}\label{defn:toric-orb}
We will denote by $\Delta^*$ an \emph{augmented fan} by which we mean a fan $\Delta$ with elements
$n(\rho) \in N\cap\rho$ for every $\rho\in\Delta(1)$.
\end{defn}
\begin{prop}\label{prop:toric-orb}
For a complete simplicial augmented fan $\Delta^*$ we have a natural orbifold structure compatible with the action of $T_\C$
on $X_\Delta$.  We denote $X_\Delta$ with this orbifold structure by $X_{\Delta^*}$.
\end{prop}
\begin{proof}
Let $\sigma\in\Delta^*(n)$ have generators $p_1,p_2,\ldots,p_n$ as in the definition.  Let
$N^\prime\subseteq N$ be the sublattice $N^\prime =\Z\{p_1,p_2,\ldots,p_n\}$, and $\sigma^\prime$ the
equivalent cone in $N^\prime$.  Denote by $M^\prime$ the dual lattice of $N^\prime$ and $T^\prime_\C$ the torus.
Then $U_{\sigma^\prime}\cong\C^n$.  It is easy to see that
\[ N/N^\prime = \Hom_\Z(M^\prime/M,\C^*). \]
And $N/N^\prime$ is the kernel of the homomorphism
\[T^\prime_\C=\Hom_\Z(M^\prime,\C^*)\rightarrow T_\C=\Hom_\Z(M,\C^*).\]
Let $\Gamma=N/N^\prime$.  An element $t\in\Gamma$ is a homomorphism $t:M^\prime\rightarrow\C^*$ equal to 1
on $M$.  The regular functions on $U_{\sigma^\prime}$ consist of $\C$-linear combinations of $x^m$ for
$m\in\sigma^{\prime\vee}\cap M^\prime$.  And $t\cdot x^m =t(m)x^m$.  Thus the invariant functions are
the $\C$-linear combinations of $x^m$ for $m\in\sigma^{\vee}\cap M$, the regular functions of $U_\sigma$.
Thus $U_{\sigma^\prime}/\Gamma = U_\sigma$.  And the charts are easily seen to be compatible on intersections.
\end{proof}
Conversely, we have the following.
\begin{prop}
Let $\Delta$ be a complete simplicial fan.  Suppose for simplicity that the local uniformizing groups are
abelian.  Then every orbifold structure on $X_\Delta$ compatible with the action of $T_\C$ arises from an
augmented fan $\Delta^*$.
\end{prop}
See~\cite{vC2} for a proof.

Let $\Delta^*$ be an augmented fan in $N$.  We will assumed from now on that the fan $\Delta$ is simplicial
and complete.
\begin{defn}
A real function $h:N_\R \rightarrow\R$ is a \emph{$\Delta^*$ -linear support function} if for
each $\sigma\in\Delta^*$ with given $\Q$-generators $p_1,\ldots,p_r$ in $N$, there is an
$l_\sigma\in M_\Q$ with $h(s)=\langle l_\sigma ,s\rangle$ and $l_\sigma$ is $\Z$-valued on the sublattice
$\Z\{p_1,\ldots,p_r\}$.  And we require that $\langle l_\sigma, s\rangle=\langle l_\tau ,s\rangle$ whenever
$s\in\sigma\cap\tau$.  The additive group of $\Delta^*$-linear support functions will be denoted by
$\SF(\Delta^*)$.
\end{defn}
Note that $h\in\SF(\Delta^*)$ is completely determined by the integers $h(n(\rho))$ for all
$\rho\in\Delta(1)$.  And conversely, an assignment of an integer to $h(n(\rho))$ for all
$\rho\in\Delta(1)$ defines $h$.  Thus
\[ \SF(\Delta^*)\cong\Z^{\Delta(1)}. \]
\begin{defn}
Let $\Delta^*$ be a complete augmented fan.  For $h\in\SF(\Delta^*)$,
\[ \Sigma_h :=\{m\in M_\R : \langle m,n\rangle \geq h(n),\text{ for all }n\in N_\R\}, \]
is a, possibly empty, convex polytope in $M_\R$.
\end{defn}

We will consider the holomorphic line $V$-bundles on $X=X_{\Delta*}$.  All $V$-bundles will be \emph{proper}
in this section.
The set of isomorphism classes of holomorphic line $V$-bundles is denoted by $\picorb(X)$, which is a group under
the tensor product.
\begin{defn}
A \emph{Baily divisor} is a $\Q$-Weil divisor $D\in\weil(X)\otimes\Q$ whose inverse image
$D_{\tilde{U}}\in\weil(\tilde{U})$ in every local uniformizing chart $\pi:\tilde{U}\rightarrow U$
is Cartier.  The additive group of Baily divisors is denoted $\divorb(X)$.
\end{defn}
A Baily divisor $D$ defines a holomorphic line $V$-bundle $[D]\in\picorb(X)$ in a way completely analogous
to Cartier divisors.  Given a nonzero meromorphic function $f\in\mathscr{M}$ we have
the \emph{principal divisor}
\[ \Div(f):=\sum_V \nu_{V}(f)V, \]
where $\nu_{V}(f)V$ is the order of the zero, or negative the order of the pole, of $f$ along
each irreducible subvariety of codimension one.
We have the exact sequence
\begin{equation}
1\rightarrow\C^* \rightarrow\mathscr{M}^*\rightarrow\divorb(X)\overset{[\ ]}{\rightarrow}\picorb(X).
\end{equation}

A holomorphic line $V$-bundle $\pi:\mathbf{L}\rightarrow X$ is equivariant
if there is an action of $T_\C$ on $\mathbf{L}$ such that $\pi$ is equivariant,
$\pi(tw)=t\pi(w)$ for $w\in\mathbf{L}$ and $t\in T_\C$ and the action lifts to a holomorphic action,
linear on the fibers, over each uniformizing neighborhood. The group of isomorphism classes of equivariant
holomorphic line $V$-bundles is denoted $\picorb_{T_\C}(X)$.
Similarly, we have invariant Baily divisors, denoted $\divorb_{T_\C}(X)$, and $[D]\in\picorb_{T_\C}(X)$
whenever $D\in\divorb_{T_\C}(X)$.
\begin{prop}\label{prop:toric-vbund}
Let $X=X_{\Delta^*}$ be compact with the standard orbifold structure, i.e. $\Delta^*$ is simplicial and complete.
\begin{thmlist}
\item  There is an isomorphism $\SF(\Delta^*)\cong\divorb_{T_\C}(X)$ obtained by sending
$h\in\SF(\Delta^*)$ to
\[ D_h :=-\sum_{\rho\in\Delta(1)}h(n(\rho))V(\rho).\]
\item  There is a natural homomorphism $\SF(\Delta^*)\rightarrow\picorb_{T_\C}(X)$ which associates \label{list:div}
an equivariant line $V$-bundle $\mathbf{L}_h$ to each $h\in\SF(\Delta^*)$.
\item  Suppose $h\in\SF(\Delta^*)$ and $m\in M$ satisfies
\[ \langle m,n\rangle\geq h(n)\text{  for all  }n\in N_\R,\]
then $m$ defines a section $\psi:X\rightarrow\mathbf{L}_h$ which has the equivariance property
$\psi(tx)=\chi^m(t)(t\psi(x))$.
\item  The set of sections $H^0(X,\mathcal{O}(\mathbf{L}_h))$ is the finite dimensional
$\C$-vector space with basis $\{x^m :m\in\Sigma_h \cap M\}$.
\item  Every Baily divisor is linearly equivalent to a $T_\C$-invariant Baily divisor.
Thus for $D\in\picorb(X)$, $[D]\cong[D_h]$ for some $h\in\SF(\Delta^*)$.
\item  If $\mathbf{L}$ is any holomorphic line $V$-bundle, then $\mathbf{L}\cong\mathbf{L}_h$
for some $h\in\SF(\Delta^*)$.  The homomorphism in part i. induces an isomorphism
$\SF(\Delta^*)\cong\picorb_{T_\C}(X)$ and we have the exact sequence
\[ 0\rightarrow M\rightarrow\SF(\Delta^*)\rightarrow\picorb(X)\rightarrow 1. \]
\end{thmlist}
\end{prop}
\begin{proof}
(i)  For each $\sigma\in\Delta(n)$ with uniformizing neighborhood
$\pi:U_{\sigma^\prime}\rightarrow U_\sigma$ as above the map
$h\rightarrow D_h$ assigns the principal divisor
\[\Div(x^{-l_\sigma})=-\sum_{\rho\in\Delta(1),\rho<\sigma}h(n(\rho))V^\prime(\rho),\]
where $V^\prime(\rho)$ is the closure of the orbit $\orb(\rho)$ in $U_{\sigma^\prime}$.
An element $\divorb_{T_\C}(X)$ must be a sum of closures of codimension one orbits $V(\rho)$ in
Proposition~\ref{prop:toric-prop}, and by above remarks the map is an isomorphism.

(ii)  One defines $\mathbf{L}_h:=[D_h]$, where $[D_h]$ is constructed as follows.  Consider a uniformizing chart
$\pi:U_{\sigma^\prime}\rightarrow U_\sigma$ as in Proposition~\ref{prop:toric-orb}.
Define $\mathbf{L}_h|_{U_{\sigma^\prime}}$ to be the invertible sheaf $\mathcal{O}_{U_{\sigma^\prime}}(D_h)$, with
$D_h$ defined on $U_{\sigma^\prime}$ by $x^{-l_\sigma}$.  So
$\mathbf{L}_h|_{U_{\sigma^\prime}}\cong U_{\sigma^\prime}\times\C$  with an action of $T^\prime_\C$,
\[ t(x,v)=(tx,\chi^{-l_\sigma}(t)v) \text{  where  }t\in T^\prime_\C,\ (x,v)\in U_{\sigma^\prime}\times\C.\]
Then $\mathbf{L}_h|_{U_\sigma}$ is the quotient by the subgroup $N/{N^\prime}\subset T^\prime_\C$, so it has
an action of $T_\C$.  And the $\mathbf{L}_h|_{U_\sigma}$ glue together equivariantly with respect to the action.

(iii)  For $\sigma\in\Delta$ we have $\langle m,n\rangle\geq\langle l_\sigma,n\rangle$ for all $n\in\sigma$.
Then $m-l_\sigma\in M^\prime \cap{\sigma^{\prime}}^\vee$ and $x^{m-l_\sigma}$ is a section of
the invertible sheaf $\mathcal{O}_{U_{\sigma^\prime}}(D_h)$ and is equivariant with respect to $N/{N^\prime}$
so it defines a section of $\mathbf{L}_h|_{U_\sigma}$. And these sections are compatible.

(iv)  We will make use of the GAGA theorems of A. Grothendieck~\cite{Gr2,Gr3}.  As with any holomorphic
$V$-bundle, the sheaf of sections $\mathcal{O}(\mathbf{L}_h)$ is a coherent sheaf.  It follows from GAGA
that we may consider $\mathcal{O}(\mathbf{L}_h)$ as a coherent algebraic sheaf, and all global sections are
algebraic.  If $\phi$ is a global section, then $\phi\in H^0(T_\C,\mathcal{O}(\mathbf{L}_h))\subset\C[M]$.
And in the uniformizing chart $\pi:U_{\sigma^\prime}\rightarrow U_\sigma$, $\phi$ lifts to an element
of the module $\mathcal{O}_{U_{\sigma^\prime}}\cdot x^{l_\sigma}$ which has a basis
$\{x^m : m\in l_\sigma +M^\prime \cap{\sigma^\prime}^\vee \}$.  So $\phi|_{U_\sigma}$ is a $\C$-linear
combination of $x^m$ with $m\in M$ and $\langle m,n\rangle\geq h(n)$ for all $n\in\sigma$. Thus
$m\in\Sigma_h$.

(v)  The divisor $T_\C \cap D$ is a Cartier divisor on $T_\C$ which is also principal since $\C[M]$ is a
unique factorization domain.  Thus there is a nonzero rational function $f$ so that
$D^\prime = D-\Div(f)$ satisfies $D^\prime \cap T_\C=\emptyset$.  Then $D^\prime \in\divorb_{T_\C}(X)$,
and the result follows from i.

(vi) Consider $\mathbf{L}_{U_{\sigma^\prime}}$ on a uniformizing neighborhood $U_{\sigma^\prime}$ as above.
For each $\rho\in\Delta(1), \rho<\sigma$ the subgroup $H_\rho\subseteq N/{N^\prime}$ fixing $V^\prime(\rho)$ is
cyclic and generated by $n^\prime \in N$ where $n^\prime$ is the primitive element with $a_\rho n^\prime =n(\rho)$.
Now $H_\rho$ acts linearly on the fibers of $\mathbf{L}_{U_{\sigma^\prime}}$ over $V^\prime(\rho)$.  Suppose
$n^\prime$ acts with weight $e^{2\pi i\frac{k}{a}}$, then let $D_\rho :=kV(\rho)$.  If
$D^\prime :=\sum_{\rho\in\Delta(1)}D_\rho$, then $\mathbf{L}^\prime :=\mathbf{L}\otimes [-D^\prime]$ is
Cartier on $X_0 :=X\setminus\Sing(X)$, where $\Sing(X)$ has codimension at least two.  The sheaf
$\mathcal{O}(\mathbf{L}^\prime)$ is not only coherent but is a rank-1 reflexive sheaf.  By GAGA
$\mathcal{O}(\mathbf{L}^\prime)\cong E\otimes\mathcal{O}^\prime$, where $E$ is an algebraic reflexive
rank-1 sheaf and $\mathcal{O}^\prime$ is the sheaf of analytic functions.  It is well known that
$E=\mathcal{O}(D)$ for $D\in\weil(X)$.  And as a Baily divisor, we have $\mathbf{L}^\prime\cong [D]$.
So $\mathbf{L}\cong [D+D^\prime]$, and by v. we have $\mathbf{L}\cong\mathbf{L}_h$ for some
$h\in\SF(\Delta^*)$.
\end{proof}

The sign convention in the proposition is adopted to make subsequent discussions involving
$\Sigma_h$ consistent with the existing literature, although having $D_{-m}=\Div(x^m)$ maybe
bothersome.  Note also that we denote a Baily divisor by a formal $\Z$-linear sum the coefficient
giving the multiplicity of the irreducible component in the \emph{uniformizing chart}.  This is different from
its expression as a Weil divisor when irreducible components are contained in codimension-1 components of the
singular set of the orbifold.

For $X=X_{\Delta^*}$ there is a unique $k\in\SF(\Delta^*)$ such that $k(n(\rho))=1$ for all $\rho\in\Delta(1)$.
The corresponding Baily divisor
\begin{equation}
D_k :=-\sum_{\rho\in\Delta(1)}V(\rho)
\end{equation}
is the \emph{(orbifold) canonical divisor}.  The corresponding $V$-bundle is $\mathbf{K}_X$, the $V$-bundle
of holomorphic n-forms.  This will in general be different from the canonical sheaf in the algebraic geometric
sense.

\begin{defn}
Consider support functions as above but which are only required to be $\Q$-valued on $N_\Q$, denoted $\SF(\Delta,\Q)$.
$h$ is \emph{strictly upper convex} if $h(n+n^\prime )\geq h(n)+h(n^\prime)$ for all $n,n^\prime \in N_\Q$
and for any two $\sigma,\sigma^\prime \in\Delta(n)$, $l_\sigma$ and $l_{\sigma^\prime}$ are different linear functions.
\end{defn}
Given a strictly upper convex support function $h$, the polytope $\Sigma_h$ is the convex hull in $M_\R$ of the
vertices $\{l_\sigma : \sigma\in\Delta(n)\}$.  Each $\rho\in\Delta(1)$ defines a facet by
\[ \langle m, n(\rho)\rangle\geq h(n(\rho)). \]
If $n(\rho)=a_\rho n^\prime$ with $n^\prime \in N$ primitive and $a_\rho\in\Z^+$ we may label the face with
$a_\rho$ to get the labeled polytope $\Sigma^*_h$ which encodes the orbifold structure.  Conversely, from
a rational convex polytope $\Sigma^*$ we associate a fan $\Delta^*$ and a support function $h$ as follows.
For an $l$-dimensional face $\theta\subset\Sigma^*$, define the rational $n$-dimensional cone
$\sigma^{\vee}(\theta)\subset M_\R$ consisting of all vectors $\lambda(p-p^\prime)$, where
$\lambda\in\R_{\geq 0}, p\in\Sigma$, and $p^\prime\in\theta$.  Then $\sigma(\theta)\subset N_\R$ is
the $(n-l)$-dimensional cone dual to $\sigma^{\vee}(\theta)$.  The set of all $\sigma(\theta)$ defines
the complete fan $\Delta^*$, where one assigns $n(\rho)$ to $\rho\in\Delta(1)$ if $n(\rho)=an^\prime$ with
$n^\prime$ primitive and $a$ is the label on the corresponding $(n-1)$-dimensional face of $\Sigma^*$.
The corresponding rational support function is then
\[h(n)=\inf\{\langle m,n\rangle : m\in\Sigma^*\}\text{  for  }n\in N_\R. \]

\begin{prop}[\cite{O2,Ful}]
There is a one-to-one correspondence between the set of pairs $(\Delta^*,h)$ with $h\in\SF(\Delta,\Q)$ strictly
upper convex, and rational convex marked polytopes $\Sigma^*_h$.
\end{prop}
We will be interested in toric orbifolds $X_{\Delta^*}$ with such a support function and polytope, $\Sigma^*_h$.
More precisely we will be concerned with the following.
\begin{defn}\label{defn:toric-Fano}
Let $X=X_{\Delta^*}$ be a compact toric orbifold.  We say that $X$ is \emph{Fano} if $-k\in\SF(\Delta^*)$,
which defines the anti-canonical $V$-bundle $\mathbf{K}^{-1}_X$, is strictly upper convex.
\end{defn}
These toric varieties are not necessarily Fano in the usual sense, since $\mathbf{K}^{-1}_X$ is the
\emph{orbifold} anti-canonical class.  This condition is equivalent to $\{n\in N_\R : k(n)\leq 1 \}\subset N_\R$
being a convex polytope with vertices $n(\rho), \rho\in\Delta(1)$.  We will use $\Delta^*$ to denote both the
augmented fan and this polytope in this case.

If $\mathbf{L}_h$ is a line $V$-bundle, then for certain $s > 0$, $\mathbf{L}_h^s\cong\mathbf{L}_{sh}$ will
be a holomorphic line bundle.  For example $s=\ord(X)$, the least common multiple of the orders of the uniformizing
groups, will do.   So suppose $\mathbf{L}_h$ is a holomorphic line bundle.  If the global holomorphic sections
generate $\mathbf{L}_h$, by Proposition
~\ref{prop:toric-vbund} $M\cap\Sigma_{h}=\{m_0,m_1,\ldots ,m_r\}$ and we have a holomorphic map
$\psi_h :X\rightarrow\cps^r$ where
\begin{equation}\label{equ:toric-map}
\psi_h(w):=[x^{m_0}(w):x^{m_1}(w):\cdots :x^{m_r}(w)].
\end{equation}
\begin{prop}[\cite{O2}]\label{prop:toric-emb}
Suppose $\mathbf{L}_h$ is a line bundle, so $h\in\SF(\Delta^*)$ is integral, and suppose $h$ is strictly
upper convex.  Then $\mathbf{L}_h$ is \emph{ample}, meaning that for large enough $\nu>0$
\[\psi_{\nu h}:X\rightarrow\cps^N ,\]
is an embedding, where $M\cap\Sigma_{\nu h}=\{m_0,m_1,\ldots,m_N\}$.
\end{prop}
\begin{cor}
Let $X$ be a Fano toric orbifold.  If $\nu>0$ is sufficiently large with $-\nu k$ integral,
$\mathbf{K}^{-\nu}$ will be very ample and $\psi_{-\nu k}:X\rightarrow\cps^N$ an embedding.
\end{cor}

Let $X_{\Delta^*}$ be an orbifold surface with $h\in\SF(\Delta^*)$.  Then the total spaces of
$\mathbf{L}_h$ and $\mathbf{L}_h^\times$ are toric varieties.  The fan of $\mathbf{L}_h^\times$ is as follows.
If $\sigma\in\Delta^*$ is spanned by $n(\rho_1),\ldots,n(\rho_k)\in\Z^n$ as in definition~\ref{defn:toric-orb} let
$\ol{\sigma}$ be the cone in $\R^{n+1}$ spanned by $(n(\rho_1),h(n(\rho_1))),\ldots,(n(\rho_k),h(n(\rho_k)))\in\Z^{n+1}$.
The collection of $\ol{\sigma}, \sigma\in\Delta^*$ defines a fan $\mathcal{C}=\mathcal{C}_h$, which is the fan of
$\mathbf{L}_h^\times$.  Furthermore, if $h$, or $-h$, is strictly upper convex, then one can add an additional $(n+1)$-cone to
$\mathcal{C}_h$ to get an affine variety $Y=\mathbf{L}_h^\times\cup\{p\}$.  We will make use of the smoothness condition on $\mathbf{L}_h^\times$.
The toric variety $\mathbf{L}_h^\times$ is smooth if for every n-cone $\ol{\sigma}$ of $\mathcal{C}$ as above
spanned by $\tau_1,\ldots,\tau_n\in\Z^{n+1}$ we have
\begin{equation}\label{eq:toric-smooth}
 (\R_{\geq 0}\tau_1 +\cdots +\R_{\geq 0}\tau_n)\cap\Z^{n+1} =\Z_{\geq 0}\tau_1 +\cdots +\Z_{\geq 0}\tau_n.
\end{equation}

Suppose $\mathbf{L}_h^q\cong\mathbf{K}_X$ for some $q\in\Z_{>0}$.  Then $\mathbf{K}_Y$ is trivial.  That is,
$Y$ has a Gorenstein singularity at the apex.

\subsection{K\"{a}hler structures}

We review the construction of toric K\"{a}hler metrics on toric varieties.  Any compact toric orbifold associated
to a polytope admits a K\"{a}hler metric (see~\cite{LT}).
Due to T. Delzant~\cite{De} and E. Lerman and S. Tolman~\cite{LT} in the orbifold case, the symplectic structure is uniquely
determined up to symplectomorphism by the polytope, which is the image of the moment map.
This polytope is $\Sigma^*_h$ of the previous section with $h$ generalized to be real valued.
There are infinitely many K\"{a}hler structures on a toric orbifold with fixed polytope $\Sigma^*_h$, but there
is a canonical K\"{a}hler metric obtained by reduction.
V. Guillemin gave an explicit formula~\cite{Gu1,CaDG} for this K\"{a}hler metric.
In particular, we show that every toric Fano orbifold admits a K\"{a}hler
metric $\omega\in c_1(X)$.

Let $\Sigma^*$ be a convex polytope in $M_\R \cong{\R^n}^*$ defined by the inequalities
\begin{equation}
\langle x,u_i\rangle\geq\lambda_i, \quad i=1,\ldots, d,
\end{equation}
where $u_i \in N\subset N_\R\cong\R^n$ and $\lambda_i \in\R$.  If $\Sigma^*_h$ is associated to $(\Delta^*,h)$, then
the $u_i$ and $\lambda_i$ are the set of pairs $n(\rho)$ and $h(n(\rho))$ for $\rho\in\Delta(1)$.
We allow the $\lambda_i$ to be real but require any set $u_{i_1},\ldots, u_{i_n}$ corresponding to a
vertex to form a $\Q$-basis of $N_\Q$.

Let $(e_1,\ldots, e_d)$ be the standard basis of $\R^d$ and $\beta:\R^d\rightarrow\R^n$ be the map
which takes $e_i$ to $u_i$.  Let $\mathfrak{n}$ be the kernel of $\beta$, so we have the exact sequence
\begin{equation}\label{eq:toric-exact}
0\rightarrow\mathfrak{n}\overset{\iota}{\rightarrow}\R^d \overset{\beta}{\rightarrow}\R^n \rightarrow 0,
\end{equation}
and the dual exact sequence
\begin{equation}\label{eq:toric-exact-dual}
0\rightarrow{\R^n}^*\overset{\beta^*}{\rightarrow}{\R^d}^*\overset{\iota^*}{\rightarrow}\mathfrak{n}^*\rightarrow 0.
\end{equation}
Since (\ref{eq:toric-exact}) induces an exact sequence of lattices, we have an exact sequence
\begin{equation}
1\rightarrow N\rightarrow T^d \rightarrow T^n \rightarrow 1,
\end{equation}
where the connected component of the identity of $N$ is an $(d-n)$-dimensional torus.
The standard representation of $T^d$ on $\C^d$ preserves the K\"{a}hler form
\begin{equation}
\frac{i}{2}\sum_{k=1}^d dz_k \wedge d\ol{z}_k,
\end{equation}
and is Hamiltonian with moment map
\begin{equation}
\mu(z)=\frac{1}{2}\sum_{k=1}^d |z_k|^2 e_k +c,
\end{equation}
unique up to a constant $c$.  We will set $c=\sum_{k=1}^d \lambda_k e_k$.  Restricting
to $\mathfrak{n}^*$ we get the moment map for the action of $N$ on $\C^d$
\begin{equation}
\mu_N (z)=\frac{1}{2}\sum_{k=1}^d |z_k|^2 \alpha_k +\lambda,
\end{equation}
with $\alpha_k =\iota^*e_k$ and $\lambda=\sum\lambda_k\alpha_k$.
Let $Z=\mu_N^{-1}(0)$ be the zero set.  By the exactness of (\ref{eq:toric-exact-dual})
$z\in\mu_N^{-1}(0)$ if an only if there is a $v\in{\R^n}^*$ with $\mu(z)=\beta^* v$.
Since $\beta^*$ is injective, we have a map
\begin{equation}
\nu: Z\rightarrow{\R^n}^*,
\end{equation}
where $\beta^*\nu(z)=\mu(z)$ for all $z\in Z$.  For $z\in Z$
\begin{equation}
\begin{split}
\langle\nu(z), u_i \rangle & =\langle\beta^*\nu(z),e_i\rangle \\
         & =\langle\mu(z), e_i \rangle \\
         & =\frac{1}{2}|z_i|^2 + \lambda_i, \\
\end{split}
\end{equation}
thus $\nu(z)\in\Sigma^*$.  Conversely, if $v\in\Sigma^*$, then $v=\nu(z)$ for some
$z\in Z$ and in fact a $T^d$ orbit in $Z$.  Thus $Z$ is compact.  The following is not
difficult to show.
\begin{thm}
The action of $N$ on $Z$ is locally free.  Thus the quotient
\[ X_{\Sigma^*} =Z/N \]
is a compact orbifold.  Let
\[ \pi: Z\rightarrow X\]
be the projection and
\[ \iota: Z\rightarrow\C^d\]
the inclusion.  Then $X_{\Sigma^*}$ has a canonical K\"{a}hler structure with K\"{a}hler form
$\omega$ uniquely defined by
\[ \pi^*\omega = \iota^*(\frac{i}{2}\sum_{k=1}^d dz_k\wedge d\ol{z}_k). \]
\end{thm}

We have an action of $T^n =T^d/N$ on $X_{\Sigma^*}$ which is Hamiltonian for $\omega$.
The map $\nu$ is $T^d$ invariant, and it descends to a map, which we also call $\nu$,
\begin{equation}
\nu:X_{\Sigma^*}\rightarrow {\R^n}^*,
\end{equation}
which is the moment map for this action.  The above comments show that
$\im(\nu)=\Sigma^*$.  The action $T^n$ extends to the complex torus $T^n_\C$ and one can show
that as an analytic variety and orbifold $X_{\Sigma^*}$ is the toric variety constructed from
$\Sigma^*$ in the previous section.  See~\cite{Gu2} for more details.

Let $\sigma :\C^d\rightarrow\C^d$ be the involution $\sigma(z)=\ol{z}$.  The set $Z$ is stable
under $\sigma$, and $\sigma$ descends to an involution on $X$.  We denote the fixed point
sets by $Z_r$ and $X_r$.  And we have the projection
\begin{equation}\label{eq:toric-proj}
\pi :Z_r \rightarrow X_r.
\end{equation}
We equip $Z_r$ and $X_r$ with Riemannian metrics by restricting the K\"{a}hler metrics on $\C^d$
and $X$ respectively.
\begin{prop}
The map (\ref{eq:toric-proj}) is a locally finite covering and is an isometry with respect to
these metrics
\end{prop}
Note that $Z_r$ is a subset of $\R^d$ defined by
\begin{equation}
\frac{1}{2}\sum_{k=1}^d x_k^2\alpha_k =-\lambda.
\end{equation}
Restrict to the orthant $x_k>0$ $k=1,\ldots,d$ of $\R^d$.  Let $Z_r^\prime$ be the component of $Z_r$ in this
orthant.  Under the coordinates
\begin{equation}
s_k = \frac{x_k^2}{2},\quad k=1,\ldots,d.
\end{equation}
The flat metric on $\R^d$ becomes
\begin{equation}\label{eq:toric-flat}
\frac{1}{2}\sum_{k=1}^d\frac{(ds_k)^2}{s_k}.
\end{equation}
Consider the moment map $\nu$ restricted to $Z_r^\prime$.  The above arguments show that
$\nu$ maps $Z_r^\prime$ diffeomorphically onto the interior $\Sigma^\circ$ of $\Sigma$.
In particular we have
\begin{equation}\label{eq:toric-real-mom}
\langle \nu(x), u_k\rangle=\lambda_k + s_k,\  k=1,\ldots,d,\text{  for  } x\in Z_r^\prime.
\end{equation}

Let $\mathit{l}_k: {\R^n}^*\rightarrow\R$ be the affine function
\[ \mathit{l}_k(x)=\langle x,u_k\rangle -\lambda_k,\quad k=1,\ldots, d.\]
Then by equation (\ref{eq:toric-real-mom}) we have
\begin{equation}
\mathit{l}_k \circ\nu = s_k.
\end{equation}
Thus the moment map $\nu$ pulls back the metric
\begin{equation}\label{eq:toric-met}
\frac{1}{2}\sum_{k=1}^d\frac{(d\mathit{l}_k)^2}{\mathit{l}_k},
\end{equation}
on $\Sigma^\circ$ to the metric (\ref{eq:toric-flat}) on $Z_r^\prime$.
We obtain the following.
\begin{prop}
The moment map $\nu :X_r^\prime\rightarrow\Sigma^\circ$ is an isometry when
$\Sigma^\circ$ is given the metric (\ref{eq:toric-met}).
\end{prop}

Let $W\subset X$ be the orbit of $T^n_{\C}$ isomorphic to $T^n_{\C}$.  Then by restriction
$W$ has a $T^n$-invariant K\"{a}hler form $\omega$.  Identify $T^n_{\C}=\C^n/{2\pi i\Z^n}$,
so there is an inclusion $\iota :\R^n\rightarrow T^n_{\C}$.
\begin{prop}\label{prop:toric-pot}
Let $\omega$ be a $T^n$-invariant K\"{a}hler form on $W$.  Then the action of $T^n$ is Hamiltonian if
and only if $\omega$ has a $T^n$-invariant potential function, that is, a function
$F\in C^{\infty}(\R^n)$ such that
\[ \omega =2i\partial\ol{\partial}F. \]
\end{prop}
\begin{proof}
Suppose the action is Hamiltonian.  Any $T^n$-orbit is Lagrangian, so $\omega$ restricts to zero.
The inclusion $T^n\subset T^n_{\C}$ is a homotopy equivalence.  Thus $\omega$ is exact.
Let $\gamma$ be a $T^n$-invariant 1-form with $\omega=d\gamma$.  Let $\gamma=\beta +\ol{\beta}$
where $\beta\in\Omega^{0,1}$.  Then
\[\omega=d\gamma =\partial\beta +\ol{\partial}\ol{\beta}, \]
since $\ol{\partial}\beta =\partial\ol{\beta}=0$.
Since $H^{0,k}(W)_{T^n}=0$ for $k>0$, there exists a $T^n$-invariant function $f$ with
$\beta=\ol{\partial}f$.  Then
\[\omega=\partial\ol{\partial}f+\ol{\partial}\partial\ol{f}=2i\partial\ol{\partial}\im f.\]
The converse is a standard result.
\end{proof}

Suppose the $T^n$ action on $W$ is Hamiltonian with moment map $\nu:W\rightarrow{\R^n}^*$.
Denote by $x+iy$ the coordinates given by the identification $W=\C^n/{2\pi i\Z^n}$.
\begin{prop}[\cite{Gu1}]\label{prop:toric-leg}
Up to a constant $\nu$ is the Legendre transform of $F$, i.e.
\[ \nu(x+iy)=\frac{\partial F}{\partial x} +c,\quad c\in{\R^n}^*\]
\end{prop}
\begin{proof}
By definition
\[d\nu_k =-\iota\biggl(\frac{\partial}{\partial y_k}\biggr)\omega.\]
But by Proposition~\ref{prop:toric-pot},
\[\omega = \sum_{j,k=1}^n \frac{\partial^2 F}{\partial x_j \partial x_k} dx_j\wedge dy_k, \]
so
\[ d\nu_k =-\iota\biggl(\frac{\partial}{\partial y_k}\biggr)\omega =d\biggl(\frac{\partial F}{\partial x_k}\biggr).\]
Therefore $\nu_k =\frac{\partial F}{\partial x_k} + c_k$.
\end{proof}
We can eliminate $c$ by replacing $F$ with $F-\sum_{k=1}^n c_k x_k$.

Notice that the metric (\ref{eq:toric-met}) on $\Sigma^\circ$ can be written
\begin{equation}
\sum_{j,k}\frac{\partial^2 G}{\partial y_j\partial y_k}dy_j dy_k,
\end{equation}
with
\begin{equation}\label{eq:toric-mompot}
G=\frac{1}{2}\sum_{k=1}^d\mathit{l}_k(y)\log\mathit{l}_k(y).
\end{equation}
V. Guillemin~\cite{Gu1} showed that the Legendre transform of $G$ is the inverse Legendre transform of $F$, i.e.
\begin{equation}
\frac{\partial F}{\partial x} =y\text{  and  }\frac{\partial G}{\partial y}=x.
\end{equation}
From this it follows that
\begin{equation}\label{eq:toric-leginv}
F(x)=\sum_{i=1}^{n}x_i y_i -G(y), \text{  where }y=\frac{\partial F}{\partial x}.
\end{equation}
Define
\[ \mathit{l}_\infty (x)=\sum_{i=1}^d\langle x,u_i\rangle. \]
From equations~(\ref{eq:toric-mompot}) and~(\ref{eq:toric-leginv}) it follows that $F$ has the expression
\begin{equation}\label{eq:toric-Fpot}
F=\frac{1}{2}\nu^*\left(\sum_{k=1}^d \lambda_k\log\mathit{l}_k +\mathit{l}_\infty\right),
\end{equation}
which gives us the following.
\begin{thm}[\cite{Gu1,CaDG}]\label{thm:toric-Gu}
On the open $T^n_{\C}$ orbit of $X_{\Sigma^*}$ the K\"{a}hler form $\omega$ is given by
\[i\partial\ol{\partial}\nu^*\left(\sum_{k=1}^d \lambda_k\log\mathit{l}_k +\mathit{l}_\infty\right).\]
\end{thm}

Suppose we have an embedding as in Proposition~\ref{prop:toric-emb},
\[ \psi_h :X_{\Sigma^*}\rightarrow\cps^N.\]
So $\Sigma_h$ is an integral polytope and $M\cap\Sigma_h=\{m_0,m_1,\ldots, m_N\}$.
Let $\omega_{\sm FS}$ be the Fubini-Study metric on $\cps^N$.  Note that $\psi_h^*\omega_{FS}$ is
degenerate along the singular set of $X$, so does not define a K\"{a}hler form.

Consider the restriction of $\psi_h$ to the open $T^n_{\C}$ orbit $W\subset X$.
Let $\iota =\psi_h |_W$.
It is induced by a representation
\begin{equation}
\tau: T^n_{\C}\rightarrow GL(N+1,\C),
\end{equation}
with weights $m_0,m_1,\ldots,m_N$.  If $z=x+iy\in\C^n/{2\pi i\Z^n} =T^n_{\C}$, and
$w=(w_0,\ldots,w_N)$, then
\begin{equation}
\tau(\exp z)w=(e^{\langle m_0,x+iy\rangle}w_0,\ldots,e^{\langle m_N,x+iy\rangle}w_N).
\end{equation}
Recall the Fubini-Study metric is
\begin{equation}
\omega_{FS} = i\partial\ol{\partial}\log |w|^2.
\end{equation}
Let $[w_0 :\cdots:w_N]$ be homogeneous coordinates of a point in the image of $W$, then
\begin{equation}
\iota^*\omega_{FS} =i\partial\ol{\partial}\log\biggl(\sum_{k=0}^N |w_k|^2 e^{2\langle m_k,x\rangle}\biggr).
\end{equation}
From equation (\ref{eq:toric-mompot}) we have
\[x =\frac{\partial G}{\partial y}=\frac{1}{2}\biggl(\sum_{j=1}^d u_j\log\mathit{l}_j +u\biggr),\]
where $u=\sum u_j$.  Then
\[ 2\langle m_i, x\rangle = 2\langle m_i,\frac{\partial G}{\partial y}\rangle=
\sum_{j=1}^d \langle m_i,u_j\rangle\log l_j +\langle m_i,u\rangle.\]
So setting $d_i=e^{\langle m_i,u\rangle}$, gives
\[ e^{2\langle m_i,x\rangle} =\nu^*\biggl(d_i\prod_{j=1}^d l_j^{\langle m_i, u_j\rangle}\biggr).\]
But from (\ref{eq:toric-Fpot}),
\[e^{2F}=\nu^*\biggl(e^{l_\infty}\prod_{j=1}^d l_j^{\lambda_j}\biggr).\]
Combining these,
\[e^{2\langle m_i,x\rangle}=e^{2F}\nu^*\biggl(d_i e^{-l_\infty}\prod_{j=1}^d l_j^{l_j(m_i)}\biggr).\]
Let $k_i =|w_i|^2 d_i$, then summing gives
\[ \sum_{i=1}^N |w_i|^2 e^{2\langle m_i,x\rangle} = e^{2F}\nu^*(e^{-l_\infty}Q),\]
where
\[ Q=\sum_{i=1}^N k_i\prod_{j=1}^d l_j^{l_j(m_i)}.\]
Thus we have
\begin{equation}\label{eq:toric-embmet}
\psi_h^*\omega_{FS} =\omega +i\partial\ol{\partial}\nu^*(-l_\infty +\log Q).
\end{equation}
Using that $\Sigma_h$ is integral, and $k_i\neq 0$ for $m_i$ a vertex of $\Sigma_h$,
it is not difficult to show that $Q$ is a positive function on $\Sigma_h$.
Thus equation (\ref{eq:toric-embmet}) is valid on all of $X$.
\begin{thm}\label{thm:toric-ampl}
Suppose $\mathbf{L}_h$ is very ample for some $h\in\SF(\Delta^*)$ strictly upper convex and integral,
and let $\omega$ be the canonical K\"{a}hler metric for the polytope $\Sigma_h$.  Then
\[ [\omega]=2\pi c_1(\mathbf{L})=[\psi_h^*\omega_{\sm FS}]. \]
\end{thm}
\begin{cor}\label{cor:toric-fano}
Suppose $X=X_{\Delta^*}$ is Fano.  Let $\omega$ be the canonical metric of the integral polytope
$\Sigma_{-k}^*$.  Then
\[ [\omega]=2\pi c_1(\mathbf{K}^{-1})=2\pi c_1(X).\]
Thus $c_1(X)>0$.  Conversely, if $c_1(X)>0$, then $\mathbf{K}^{-p}$ is very ample for some $p>0$
and $X$ is Fano as defined in definition (\ref{defn:toric-Fano}).
\end{cor}
\begin{proof}
For some $p\in\Z^+$, $-pk\in\SF(\Delta^*)$ is integral and
$\mathbf{L}_{-pk}=\mathbf{K}^{-p}$ is very ample.  Let $\tilde{\omega}$ be the canonical metric of the integral polytope
$\Sigma_{-pk}^*$.  From the theorem we have
\[ [\tilde{\omega}]=2\pi c_1(\mathbf{K}^{-p})=2\pi pc_1(X).\]
Let $\omega$ be the canonical metric for $\Sigma^*_{-k}$.
Theorem~(\ref{thm:toric-Gu}) implies that $[\tilde{\omega}]=p[\omega]$

For the converse, It follows from the
extension to orbifolds of the Kodaira embedding theorem of W. Baily~\cite{Ba} that
$\mathbf{K}^{-p}$ is very ample for some $p>0$ sufficiently large.  It follows from
standard results on toric varieties that $-k$ is strictly upper convex (see~\cite{O2}).
\end{proof}

The next result will have interesting applications to the Einstein manifolds constructed later.
\begin{prop}
With the canonical metric the volume of $X_{\Sigma^*}$ is $(2\pi)^n$ times the Euclidean volume
of $\Sigma$.
\end{prop}
\begin{proof}
Let $W\subset X$ be the open $T_\C^n$ orbit.
We identify $W$ with $\C^n/{2\pi i\Z^n}$ with coordinates $x+iy$.
The restriction of $\omega$ to $W$ is
\[\omega|_W = \sum_{j,k=1}^n \frac{\partial^2 F}{\partial x_j \partial x_k} dx_j\wedge dy_k. \]
Thus
\[\frac{\omega^n}{n!}=\det\biggl(\frac{\partial^2 F}{\partial x_j \partial x_k}\biggr)dx\wedge dy.\]
Integrating over $dy$ gives
\[ \Vol(X,\omega)=(2\pi)^n \int_{\R^n}\det\biggl(\frac{\partial^2 F}{\partial x_j \partial x_k}\biggr)dx.\]
By Proposition~\ref{prop:toric-leg} $x\rightarrow z=\nu(x+iy)=\frac{\partial F}{\partial x}$ is a diffeomorphism
from $\R^n$ to $\Sigma^\circ$.  By the change of variables,
\[\Vol(\Sigma)=\int_{\Sigma}dz=\int_{\R^n}\det\biggl(\frac{\partial^2 F}{\partial x_j \partial x_k}\biggr)dx.\]
\end{proof}
\begin{cor}\label{cor:toric-fano-vol}
Let $X=X_{\Delta^*}$ be a toric Fano orbifold.  And let $\omega$ be any K\"{a}hler form with
$\omega\in c_1(X)$.  Then
\[ \Vol(X,\omega)=\frac{1}{n!}c_1(X)^n[X]=\Vol(\Sigma_{-k}).\]
\end{cor}
\begin{proof}
Let $\omega_{\sm c}$ be the canonical metric associated to $\Sigma_{-k}^*$,
then $\frac{1}{2\pi}\omega_{\sm c} \in c_1(X)$ by Corollary~\ref{cor:toric-fano}.
Then
\[\Vol(X,\omega)=\frac{1}{(2\pi)^n}\Vol(X,\omega_{\sm c})=\Vol(\Sigma_{-k}).\]
\end{proof}

\subsection{Moment map and Futaki invariant}

A closer analysis of the moment map in the Fano case, originally due to T. Mabuchi~\cite{Mab},
will be useful.  Suppose in this section that $X=X_{\Delta^*}$ is toric Fano with $\dim_\C X=n$.
As above,  $z_k=x_k +iy_k, 1\leq k\leq n$ are logarithmic coordinates,
$(z_1,\ldots,z_n)\in\C^n/{2\pi i\Z^n}\cong T_\C$.
For $(t_1,\ldots,t_n)\in T_\C$, $x_k =\log|t_k|$.  Then a $T$-invariant function $u\in C^\infty(T_\C)$ is considered as
a $C^\infty$ function $u=u(x_1,\ldots,x_n)$ on $\R^n$.  There exists a $T$-invariant fiber metric
$\Omega$ on $\mathbf{K}_X^{-1}$ with positive Chern form.  Thus there exits a $C^\infty$ function $u=u(x_1,\ldots,x_n)$
so that
\begin{equation}
e^{-u}\prod_{k=1}^n (dx_k\wedge dy_k),
\end{equation}
extends to a volume form $\Omega$ on all of $X$ and $i\partial\ol{\partial}u$ extends to a K\"{a}hler form $\omega$.
The moment map $\nu_u:X\rightarrow M_\R$ can be given, without an ambiguous constant, as
\begin{equation}
\nu_u(t) = (\frac{\partial u}{\partial x_1}(t),\ldots,\frac{\partial u}{\partial x_n}(t)),\quad\text{for }t\in T_\C.
\end{equation}
\begin{thm}\label{thm:moment-map}
The closure of the image $\nu_u(T_\C)$ in $M_\R$ is $\Sigma_{-k}$.  Furthermore, $\nu_u$ extends to a $C^\infty$
map $\nu_u: X\rightarrow M_\R$, which is the usual moment map.
\end{thm}
Only the first statement remains to be proved.  This is a slight generalization of a similar result in~\cite{Mab}.

We define the Futaki invariant.   Let $\Aut^o(X)\subseteq\Aut(X)$ be the subgroup of the homomorphic automorphism
group preserving the orbifold structure, and let $\mathfrak{g}$ be its Lie algebra.
Let $\omega\in 2\pi c_1(X)$ be a K\"{a}hler form.  There exists $f\in C^\infty (X)$ with
$\ric(\omega)-\omega =i\partial\ol{\partial}f$.  Set $c=-2^{n+1}((2\pi c_1(X))^n[X])^{-1}$.  Then the
\emph{Futaki invariant} $F:\mathfrak{g}\rightarrow\C$ is defined by
\begin{equation}\label{Fut-inv}
 F(V)=c\int_X Vf\omega^n, \quad\text{for }V\in\frak{g}.
\end{equation}
Note that, as proved in~\cite{Fut}, $F$ is zero on $[\frak{g},\frak{g}]$.  We have the Cartan decomposition
\[\frak{g}=\frak{t}_\C \oplus\sum_i\C v_i,\]
where the $v_i$ are eigenvalues for the adjoint action of $\frak{t}_\C$.  Since $\frak{t}$ is the Lie algebra of a
maximal torus, one see that the $v_i$ are contained in $[\frak{g},\frak{g}]$.  Thus we may restrict $F$ to
$\frak{t}_\C$.

Suppose $\mathbf{L}$ is an equivariant holomorphic line V-bundle on $X$.  In our case of interest $\mathbf{L}=\mathbf{K}_X^{-1}$
with the usual action.  Let $H$ be the space of hermitian metrics on $\mathbf{L}$.  For $h\in H$ denote by
$c_1(\mathbf{L},h)$ the Chern form, $\frac{i}{2\pi}\ol{\partial}\partial\log(h)$ in local holomorphic coordinates.
For a pair $(h',h'')\in H\times H$ we define
\[ R_{\mathbf{L}}(h',h'') :=\int_{t_0}^{t_1}\left(\int_X h_t^{-1}\dot{h}_t(2\pi c_1(\mathbf{L},h_t))^n \right)dt, \]
where $h_t, t_0\leq t_1$ is any piecewise smooth path with $h_{t_0}=h'$ and $h_{t_1}=h''$.
One has that $R_{\mathbf{L}}(h',h'')$ is independent of the path $h_t$ and satisfies
\[ R_{\mathbf{L}}(g^*h',g^*h'')=R_{\mathbf{L}}(h',h''),\quad\text{for }g\in\Aut^o(X),\]
and the cocycle conditions
\begin{gather}
R_{\mathbf{L}}(h',h'')+R_{\mathbf{L}}(h'',h')=0, \text{and} \\
R_{\mathbf{L}}(h,h')+R_{\mathbf{L}}(h',h'')+R_{\mathbf{L}}(h'',h)=0,
\end{gather}
for any $h,h',h''\in H$.  These identities imply that
\begin{equation}\label{Futak:char}
r_{\mathbf{L}}(g) := \exp(R_{\mathbf{L}}(h,g^*h),\quad\text{for }g\in\Aut^o(X),
\end{equation}
is in independent of $h\in H$ and is a Lie group homomorphism into $\R_+$.  It has associated Lie algebra
character ${r_{\mathbf{L}}}_* :\frak{g}\rightarrow\R$.

Let $\sigma$ be a $T_\C$-invariant section of $\mathbf{L}^*$.  Then $h\in H$ is
$h=e^{-u_h}\sigma\otimes\ol{\sigma}$ on $W\cong T_\C$ for some $u_h\in C^\infty(W)$.  We denote
by $V_\R$ the real component of a homomorphic vector field.  Then differentiating (\ref{Futak:char}) gives
\begin{equation}\label{Futak:diff}
{r_{\mathbf{L}}}_*(V) =-\int_W V_{\R}(u_h)(i\partial\ol{\partial}u_h)^n,
\end{equation}
for $V\in\frak{t}_\C$ independent of $h\in H$.  We have the following.
\begin{prop}
Suppose $X_{\Delta^*}$ is a toric Fano orbifold.  Then with $\mathbf{L}=\mathbf{K}_X^{-1}$ we have
\[ F_X = -2^{n+1}((2\pi c_1(X))^n[X])^{-1} {r_{\mathbf{L}}}_*,\]
where both sides are restricted to $\frak{t}_\C$.
\end{prop}
\begin{proof}
 By assumption there is an $h\in H$ with positive Chern form.  Let $\beta$ be a $T_\C$-invariant section
 of $\mathbf{K}_X^{-1}$ over $W$.  Then $h$ may be written as the volume form on $W$
 $\Omega =i^n(-1)^{n(n-1)/2}e^{-u_h}\beta\wedge\ol{\beta}$ and $i\partial\ol{\partial}u_h$ extends to a
 K\"{a}hler form $\omega\in 2\pi c_1(X)$.  Then if $f=\log(\frac{\Omega}{\omega^n})$, we have
 \[\ric(\omega)-\omega =i\partial\ol{\partial}f.\]
 Then
 \[\begin{split}
  0 & =\int_X \mathscr{L}_{V_\R}\left( e^{-f}\Omega\right) \\
    & =-\int_X V_\R(f)\omega^n +\int_X e^{-f}\mathscr{L}_{V_\R}\Omega \\
    & =-\int_X V_\R(f)\omega^n -\int_W V_{\R}(u_h)\omega^n
 \end{split}\]
And the result follows from (\ref{Futak:diff}).
\end{proof}

The $t_k\frac{\partial}{\partial t_k}, k=1,\ldots,n$ from a basis of $\frak{t}_\C$.  We may assume that the K\"{a}hler
form $\omega$ is $T$-invariant.  Then we have
\[\begin{split}
   F_X(t_k\frac{\partial}{\partial t_k}) & =-2^{n+1}((2\pi c_1(X))^n[X])^{-1} {r_\mathbf{L}}_*(t_k\frac{\partial}{\partial t_k})\\
    & =2^n((2\pi c_1(X))^n[X])^{-1}\int_W \frac{\partial u_h}{\partial x_k}(i\partial\ol{\partial}u_h)^n\\
    & =2^n((2\pi c_1(X))^n[X])^{-1}\int_W \frac{\partial u_h}{\partial x_k}\left(\frac{1}{2}\right)^n n!
    \det\left(\frac{\partial^2 u_h}{\partial x_j \partial x_k}\right)\prod_{l=1}^n dx_l \wedge dy_l \\
    & =(2\pi)^n((2\pi c_1(X))^n[X])^{-1}\int_{\R^n}\frac{\partial u_h}{\partial x_k}n!
    \det\left(\frac{\partial^2 u_h}{\partial x_j \partial x_k}\right)dx\\
    & =(2\pi)^n n!((2\pi c_1(X))^n[X])^{-1}\int_{\Sigma_{-k}} y_k dy\\
    & =\Vol(\Sigma_{-k})^{-1}\int_{\Sigma_{-k}} y_k dy\\
  \end{split}\]
where $y_k=\frac{\partial u_h}{\partial x_k}$.

We have the following simple interpretation of the Futaki invariant in this case.
\begin{prop}\label{prop:Fut-bary}
Suppose $X=X_{\Delta^*}$ is a toric Fano orbifold.  Then the Futaki invariant $F_X$ is the barycenter
of the polytope $\Sigma_{-k}$.
\end{prop}

\subsection{Symmetric toric orbifolds}

Let $X_\Delta$ be an $n$-dimensional toric variety.
Let $\mathcal{N}(T_{\C})\subset\Aut(X)$ be the normalizer of $T_{\C}$.  Then
$\mathcal{W}(X):=\mathcal{N}(T_{\C})/T_{\C}$ is isomorphic to the finite group of all symmetries of $\Delta$,
i.e. the subgroup of $GL(n,\Z)$ of all $\gamma\in GL(n,\Z)$ with $\gamma(\Delta)=\Delta$.
Then we have the exact sequence.
\begin{equation}\label{toric:wely-exact1}
1\rightarrow T_{\C}\rightarrow\mathcal{N}(T_{\C}) \rightarrow\mathcal{W}(X)\rightarrow 1.
\end{equation}
Choosing a point $x\in X$ in the open orbit, defines an inclusion $T_{\C}\subset X$.  This also
provides a splitting of (\ref{toric:wely-exact1}).
Let $\mathcal{W}_0(X)\subseteq\mathcal{W}(X)$ be the subgroup which are also automorphisms of $\Delta^*$;
$\gamma\in\mathcal{W}_0(X)$ is an element of $\mathcal{N}(T_{\C})\subset\Aut(X)$ which preserves the orbifold
structure.  Let $G\subset\mathcal{N}(T_{\C})$ be the compact subgroup generated by $T^n$, the maximal compact subgroup
of $T_\C$, and $\mathcal{W}_0(X)$.
Then we have the, split, exact sequence
\begin{equation}\label{weyl:exact2}
1\rightarrow T^n\rightarrow G\rightarrow\mathcal{W}_0(X)\rightarrow 1.
\end{equation}
\begin{defn}
A \emph{symmetric Fano toric orbifold} $X$ is a Fano toric orbifold with
$\mathcal{W}_0$ acting on $N$ with the origin as the only fixed point.
Such a variety and its orbifold structure is characterized by the convex polytope $\Delta^*$ invariant under $\mathcal{W}_0$.
We call a toric orbifold \emph{special symmetric} if $\mathcal{W}_0(X)$ contains the
involution $\sigma: N\rightarrow N$, where $\sigma(n)=-n$.
\end{defn}
The following is immediate from Proposition~\ref{prop:Fut-bary}.
\begin{prop}
 For a symmetric Fano toric orbifold $X$ one has $F_X\equiv 0$.
\end{prop}

\begin{defn}\label{defn:toric-ind}
The \emph{index} of a Fano orbifold $X$ is the largest positive integer $m$ such that there is a holomorphic
$V$-bundle $\mathbf{L}$ with $\mathbf{L}^m \cong\mathbf{K}_X^{-1}$. The index of $X$ is denoted
$\ind(X)$.
\end{defn}
Note that $c_1(X)\in H^2_{orb}(X,\Z)$, and $\ind(X)$ is the greatest positive integer $m$ such that
$\frac{1}{m}c_1(X)\in H^2_{orb}(X,\Z)$.
\begin{prop}\label{prop:toric-ind-sp}
Let $X_{\Delta^*}$ be a special symmetric toric Fano orbifold.  Then $\ind(X)=1$ or $2$.
\end{prop}
\begin{proof}
We have $\mathbf{K}^{-1}\cong\mathbf{L}_{-k}$ with $-k\in\SF(\Delta^*)$ where $-k(n_\rho)=-1$
for all $\rho\in\Delta(1)$.  Suppose we have $\mathbf{L}^m\cong\mathbf{K}^{-1}$.
By Proposition~\ref{prop:toric-vbund} there is an $h\in\SF(\Delta^*)$ and $f\in M$ so
that $mh=-k+f$.  For some $\rho\in\Delta(1)$,
\begin{gather*}
mh(n_\rho) = -1 + f(n_\rho) \\
mh(-n_\rho) = -1 - f(n_\rho).
\end{gather*}
Thus $m(h(n_\rho)+h(-n_\rho))=-2$, and $m=1$ or $2$.
\end{proof}

\begin{figure}
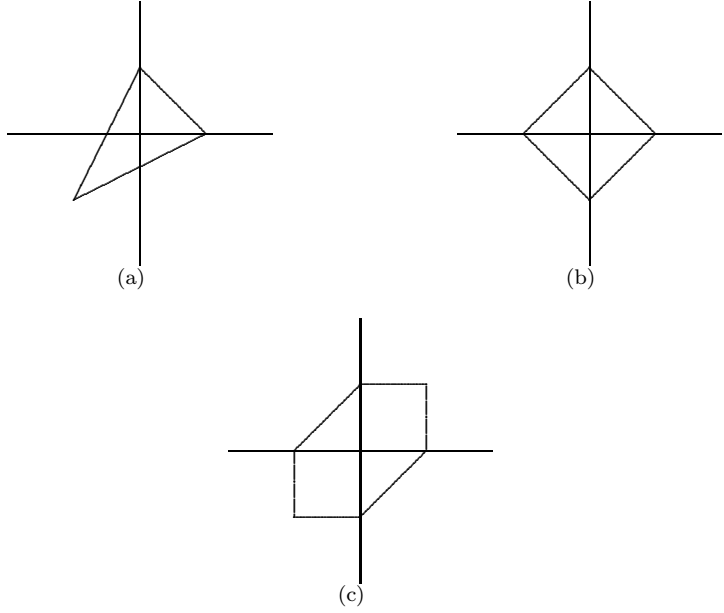

\centering

\mbox{
\subfigure[]{
\beginpicture
\setcoordinatesystem units <.5pt,.5pt> point at 0 0
\setplotarea x from -100 to 100, y from -100 to 100
\axis bottom shiftedto x=0 /
\axis left shiftedto y=0 /
\setlinear
\plot 0 50  50 0  -50 -50   0 50 /
\endpicture}

\hspace{60pt}

\subfigure[]{
\beginpicture
\setcoordinatesystem units <.5pt,.5pt> point at 0 0
\setplotarea x from -100 to 100, y from -100 to 100
\axis bottom shiftedto x=0 /
\axis left shiftedto y=0 /
\setlinear
\plot 0 50   50 0  0 -50  -50 0  0 50 /
\endpicture}}

\subfigure[]{
\beginpicture
\setcoordinatesystem units <.5pt,.5pt> point at 0 0
\setplotarea x from -100 to 100, y from -100 to 100
\axis bottom shiftedto x=0 /
\axis left shiftedto y=0 /
\setlinear
\plot 0 50  50 50  50 0  0 -50  -50 -50  -50 0  0 50 /
\endpicture}

\caption{The three smooth examples}\label{fig:sm-xpl}
\end{figure}

We will now restrict to dimension two.
In the smooth case every Fano surface, called a \emph{del Pezzo surface}, is either $\cps^1\times\cps^1$ or
$\cps^2$ blown up at $r$ points in general position $0\leq r\leq 8$.  The smooth toric Fano surfaces are
$\cps^1\times\cps^1$, $\cps^2$, the Hirzebruch surface $F_1$, the equivariant blow up of $\cps^2$ at
two $T_{\C}$-fixed points, and the equivariant blow up of $\cps^2$ at three $T_{\C}$-fixed points.
There are only three examples of smooth symmetric toric Fano surfaces, which are $\cps^1\times\cps^1$,
$\cps^2$, and the equivariant blow up of $\cps^2$ at three $T_{\C}$-fixed points.  Their marked fans
are shown in figure~\ref{fig:sm-xpl}.
The smooth toric Fano surfaces admitting a K\"{a}hler-Einstein metric are precisely the symmetric cases.

\begin{figure}
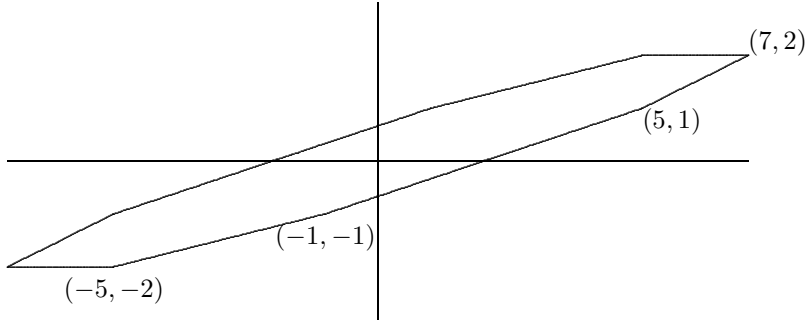

\centering
\mbox{
\beginpicture
\setcoordinatesystem units <.4pt,.4pt> point at 0 0
\setplotarea x from -350 to 350, y from -150 to 150
\axis bottom shiftedto x=0 /
\axis left shiftedto y=0 /
\setlinear
\plot 50 50  250 100  350 100  250 50  -50 -50  -250 -100  -350 -100  -250 -50  50 50 /
\put {$(7,2)$} [lb] at 350 100
\put {$(5,1)$} [lt] at 250 50
\put {$(-1,-1)$} [ct] at -50 -60
\put {$(-5,-2)$} [ct] at -250 -110
\endpicture}

\caption{Example with 8 point singular set and $\mathcal{W}_0=\Z_2$}\label{fig:typ-xpl}
\end{figure}

\section{Einstein equation}\label{sec:Einst-eq}

We consider the existence of Sasaki-Einstein metrics on toric Sasakian manifolds.  This problem is completely solved in
~\cite{FOW} where the more generally the existence of Sasaki-Ricci solitons is proved extending the existence of K\"{a}hler-Ricci
solitons on toric Fano manifolds proved in~\cite{WaZhu}.
For the examples of Sasaki-Einstein manifolds considered in this article a complete proof can be found in~\cite{vC2}.

Given a Sasakian manifold $(M,g,\Phi,\xi,\eta)$ we consider deformations of the transversal K\"{a}hler structure.  That is,
for a basic function $\phi\in C^\infty_{B}(M)$, set
\[\tilde{\eta} =\eta +2d^c_B \phi.\]
Then for $\phi$ small enough,
\[\tilde{\omega}^T =\frac{1}{2}d\tilde{\eta}=\frac{1}{2}d\eta +d_{B}d^c_{B}\phi =\omega^T +d_{B}d^c_{B}\phi \]
is a transversal K\"{a}hler metric and $\tilde{\eta}\wedge\tilde{\eta}^m$ is nowhere zero.  Then there is a Sasakian
structure $(M,\tilde{g},\tilde{\Phi},\xi,\tilde{\eta})$ with the same Reeb vector field, transverse holomorphic structure,
and basic K\"{a}hler class $[\tilde{\omega}^T]_B=[\omega^T]_B$.
The existence of such a deformation of the transverse K\"{a}hler structure to a Sasakian Einstein structure requires
the following.
\begin{prop}\label{CY-cond}
The following three conditions are equivalent.
\begin{thmlist}
\item  $(2m+2)\omega \in 2\pi c_1^B(M).$

\item  $\mathbf{K}_{C(M)}^q$ is holomorphically trivial for some $q\in\Z_+$, and there is a nowhere zero section
$\Omega$ of $\mathbf{K}_{C(M)}$, which will be multi-valued if $\mathbf{K}_{C(M)}^q$ is not trivial,
with $\mathscr{L}_{r\partial_r} \Omega =(m+1)\Omega$.

\item  $\mathbf{K}_{C(M)}^q$ is holomorphically trivial for some $q\in\Z_+$, and there is a section $\Omega$ of
$\mathbf{K}_{C(M)}$, which will be multi-valued if $\mathbf{K}_{C(M)}^q$ is not trivial, such that
\[\frac{i^{m+1}}{2^{m+1}}(-1)^{\frac{m(m+1)}{2}}\Omega\wedge\ol{\Omega} =e^f \frac{1}{n!}\ol{\omega}^{m+1},\]
for $f\in C^\infty_B(M)$ pulled back to an element of $C^\infty(C(M))$.  Here $\ol{\omega}$ is the K\"{a}hler form on $C(M)$.
\end{thmlist}
\end{prop}
In Proposition~\ref{KE-Sasak} with $\mathbf{L}=\mathbf{K}_X^{\frac{1}{q}}$ the conditions of Proposition~\ref{CY-cond}
are satisfied.

We will need the following definition.
\begin{defn}
A \emph{Hamiltonian holomorphic vector field} on $M$  is a complex vector field $Y$ invariant by $\xi$ so that
\begin{thmlist}
 \item  For any local leaf space projection $\pi_\alpha :U_\alpha \rightarrow V_\alpha$, $\pi_\alpha(Y)$ is a holomorphic
 vector field,
 \item  the complex function $\theta_Y =\sqrt{-1}\eta(Y)$ satisfies
 \[ \ol{\partial}_B \theta_Y =-\frac{\sqrt{-1}}{2}Y\neg d\eta.\]
\end{thmlist}
The Lie algebra of Hamiltonian holomorphic vector fields on $M$ is denoted $\mathfrak{h}$.
\end{defn}
Note that $\tilde{Y}=Y+i\eta(Y)r\partial_r$ is a holomorphic vector field on $C(M)$.  These correspond exactly to
transversely holomorphic vector fields, i.e. those satisfying (i), which have a potential function (cf.~\cite{BGS}).
In other words if $Y\in\Gamma(D\otimes\C)$ satisfies (i) and $\mathscr{L}_\xi Y=0$, and there exits
a complex function $\theta_Y\in C^\infty_B(M)$ with
\[ \ol{\partial}_B \theta_Y =-\frac{\sqrt{-1}}{2}Y\neg d\eta,\]
then $Y-\sqrt{-1}\theta_Y \xi$ is a Hamiltonian holomorphic vector field.
Furthermore, in the case $c_1^B(M)>0$ by the transverse Calabi-Yau theorem (cf.~\cite{EKN})
there exits a transversal K\"{a}hler deformation
$(M,\tilde{g},\tilde{\Phi},\xi,\tilde{\eta})$ with $\tilde{\eta} =\eta +2d^c_B \phi$ for some basic $\phi\in C^\infty_B(M)$
with $\Ric^T_{\tilde{g}^T}$ positive.  The usual Weitzenb\"{o}ck formula shows that the space of basic harmonic 1-forms
$\mathscr{H}^1_B$ is zero.  It follows that if $Y\in\Gamma(D\otimes\C)$ satisfies (i) and $\mathscr{L}_\xi Y=0$, there exits
a potential function $\theta_Y\in C^\infty_B(M)$ so that $Y-\sqrt{-1}\theta_Y \xi$
is a Hamiltonian holomorphic vector field.  Thus $\mathfrak{h}$ is isomorphic to the space of transversely holomorphic
vector fields commuting with $\xi$, which is isomorphic to the space of holomorphic vector fields on $C(M)$
commuting with $\xi +\sqrt{-1}r\partial_r$.

We now suppose that the conditions of Proposition~\ref{CY-cond} hold for $(M,g,\Phi,\xi,\eta)$.
Let $Y$ be a Hamiltonian holomorphic vector field.  Then $(M,g,\Phi,\xi,\eta)$ is a
\emph{Sasakian-Ricci soliton} if
\begin{equation}\label{Sasak-Ric-sol}
 \Ric^T - (2m+2)g^T =\mathscr{L}_Y g^T.
\end{equation}

Let $h\in C^\infty_B(M)$ be a basic function with
\begin{equation}
 \ric(\omega^T) -(2m+2)\omega^T =i\partial_B \ol{\partial}_B h,
\end{equation}
where $\omega^T =\frac{1}{2}d\eta$ is the transverse K\"{a}hler form.
In~\cite{TZh} Tian and Zhu defined a modified Futaki invariant $F_Y$
\begin{equation}\label{eq:Futak-modified}
 F_Y(v) =\int_M v(h -\theta_Y)e^{\theta_Y}\eta\wedge(\frac{1}{2}d\eta)^m ,\quad v\in\frak{h}.
\end{equation}
One can show as in~\cite{TZh} that $F_Y$ is unchanged under transversal deformations
$\eta \rightarrow \tilde{\eta} +2d^c_B \phi$ of the Sasakian structure; and
$F_Y(v) =0$, for all $v\in\mathfrak{h}$, is a necessary condition for a solution to (\ref{Sasak-Ric-sol}).
If $\theta_Y$ is a constant, i.e. $Y=c\xi$ for some $c\in\C$, then (\ref{eq:Futak-modified}) defines the usual Futaki
invariant.  And if $M$ is quasi-regular this is, up to a constant, the same invariant defined in (\ref{Fut-inv}).

Let $H$ be the transversal holomorphic automorphism group generated by $\mathfrak{h}$.  Let $K\subset H$ be
a compact group with Lie algebra $\mathfrak{k}$, and let $\mathfrak{k}^\C$ be its complexification.
Note that we may choose a $K$-invariant Sasakian structure $(M,\tilde{g},\tilde{\Phi},\xi,\tilde{\eta})$.  Then
for $Y\in\mathfrak{k}$, one can take $\theta_Y$ to be imaginary.  So one has $\mathfrak{k}^\C\subset\mathfrak{h}$.
As in~\cite{TZh} we have the following.
\begin{prop}
 There exists a $Y\in\mathfrak{k}^\C$ with $\im Y\in\mathfrak{k}$ so that
 \[ F_Y(v) =0, \quad\text{for all } v\in\mathfrak{k}^\C .\]
 Furthermore, $Y$ is unique up to addition of $c\xi$, for $c\in\C$, and
 \[ F_Y([v,w])=0, \quad\text{for all }v\in\mathfrak{k}^\C\text{ and }w\in\mathfrak{h}.\]
\end{prop}
Suppose now that $(M,g,\Phi,\xi,\eta)$ is a toric Sasakian manifold.  If $\mathfrak{t}$ is the Lie algebra of the
$m+1$-torus $T$ acting on $M$, then $\mathfrak{t}\subset\mathfrak{t}_\C\subseteq\mathfrak{h}$.
Using the same argument as after (\ref{Fut-inv}) we have the following.
\begin{cor}\label{cor:uni-vf}
 If $M$ is toric, then there exists a unique $Y\in\mathfrak{t}_\C$ with $\im Y\in\mathfrak{t}$ so that
 \[ F_Y(v) =0, \quad\text{for all }v\in\mathfrak{h}.\]
\end{cor}

For $\phi\in C^\infty_{B}(M)$, $\tilde{\eta} =\eta +2d^c_B \phi$ defines a transversally deformed Sasakian structure,
with transverse K\"{a}hler form $\tilde{\omega}^T =\omega +i\partial\ol{\partial}\phi$. And the Hamiltonian function for
$Y\in\mathfrak{h}$ becomes $\tilde{\theta}_Y =\theta_Y +Y\phi$ (see~\cite{TZh}).  In transverse holomorphic
coordinates (\ref{Sasak-Ric-sol}) becomes the Monge-Amp\`{e}re equation
\begin{equation}\label{Monge-Amp}
 \frac{\det(g^T_{i\ol{j}} + \phi_{i\ol{j}})}{\det(g^T_{i\ol{j}})}=e^{-(2m+2)\phi -\theta_Y -Y\phi +h}.
\end{equation}
In the toric case A. Futaki, H. Ono, and G. Wang~\cite{FOW} prove the necessary $C^0$ estimate on $\phi$ to solve
(\ref{Monge-Amp}) using the continuity method.
\begin{thm}\label{thm:Sasak-Ric-sol}
Let $(M,g,\Phi,\xi,\eta)$ be a compact toric Sasakian manifold satisfying the conditions in (\ref{CY-cond}).
Then there exists a unique transversal deformation $(M,\tilde{g},\tilde{\Phi},\xi,\tilde{\eta})$ which is a Sasakian-Ricci soliton.
\end{thm}
If the Futaki invariant vanishes, then one has $Y=0$ in Corollary~\ref{cor:uni-vf}.  Therefore we have the following.
\begin{cor}\label{cor:Sasak-Einst}
The solution in Theorem~\ref{thm:Sasak-Ric-sol} is Einstein if, and only if, the Futaki invariant vanishes.
\end{cor}

If $(M,g,\Phi,\xi,\eta)$ is quasi-regular then the leaf space of $\mathscr{F}_\xi$ is a
toric orbifold $X=X_{\Sigma^*_{-k}}$ where $\Sigma^*_{-k}$ is the marked polytope of the anti-canonical bundle.
Then $(M,g,\Phi,\xi,\eta)$ admits a transversal deformation to an Sasaki-Einstein structure if, and only if,
the barycenter of $\Sigma^*_{-k}$ is the origin.  In particular, if $M$ is quasi-regular and the leaf space
$X_{\Delta^*}$ is symmetric, then solution is Einstein.

\section{3-Sasakian manifolds}\label{sec:3Sasak}

In this section we define 3-Sasakian manifolds, the closely related quaternionic-K\"{a}hler spaces, and their
twistor spaces.  For more details see~\cite{BG2}.
These are sister geometries where one is able to pass from one to the other two by considering
the appropriate orbifold fibration.  Given a 3-Sasakian manifold $\mathcal{S}$ there is the associated
twistor space $\mathcal{Z}$, quaternionic-K\"{a}hler orbifold $\mathcal{M}$, and hyperk\"{a}hler
cone $C(\mathcal{S})$.
This is characterized by the \emph{diamond}:

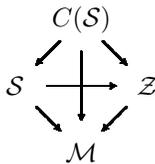
\begin{center}
\mbox{
\beginpicture
\setcoordinatesystem units <1pt, 1pt> point at 0 0
\put {$C(\mathcal{S})$} at 0 25
\put {$\mathcal{S}$} at -25 0
\put {$\mathcal{Z}$} at 25 0
\put {$\mathcal{M}$} at 0 -25
\plotsymbolspacing=.01pt
\arrow <2pt> [.3, 1] from 8 17 to 17 8
\arrow <2pt> [.3, 1] from -8 17 to -17 8
\arrow <2pt> [.3, 1] from -17 -8 to -8 -17
\arrow <2pt> [.3, 1] from 17 -8 to 8 -17
\arrow <2pt> [.3, 1] from -13 0 to 13 0
\arrow <2pt> [.3, 1] from 0 13 to 0 -13
\endpicture}
\end{center}

The equivalent 3-Sasakian and quaternionic-K\"{a}hler
reduction procedures provide an elementary method for constructing 3-Sasakian and quaternionic-K\"{a}hler
orbifolds (cf.~\cite{GL,BGMR}).  This method is effective in producing smooth 3-Sasakian manifolds, though the quaternionic-K\"{a}hler
spaces obtained are rarely smooth.  In particular, we are interested in toric 3-Sasakian 7-manifolds $\mathcal{S}$ and their
associated four dimensional quaternionic-K\"{a}hler orbifolds $\mathcal{M}$.  Here toric means that the structure is preserved
by an action of the real two torus $T^2$.  In four dimensions quaternionic-K\"{a}hler means that $\mathcal{M}$ is
Einstein and anti-self-dual, i.e. the self-dual half of the Weyl curvature vanishes $W_+\equiv 0$.
These examples are well known and they are all obtained by reduction (cf. ~\cite{BGMR} and~\cite{CaSin}).
In this case we will associate two more Einstein spaces to the four Einstein spaces in the diamond.
To each diamond of a toric 3-Sasakian manifold we have a special symmetric toric Fano surface $X$
and a Sasaki-Einstein manifold $M$ which complete diagram~\ref{int:diag-cor}.

The motivation is two fold.  First, it adds two more Einstein spaces to the examples on the right considered by
C. Boyer, K. Galicki, and others in~\cite{BGMR,BG2} and also by D. Calderbank and M. Singer~\cite{CaSin}.
Second, $M$ is smooth when the 3-Sasakian space $\mathcal{S}$ is.  And the smoothness of $\mathcal{S}$ is ensured by a relatively
mild condition on the moment map.  Thus we get infinitely many quasi-regular Sasaki-Einstein 5-manifolds with arbitrarily
high second Betti numbers paralleling the 3-Sasakian manifolds constructed in~\cite{BGMR}.

\subsection{Definitions and basic properties}\label{subsec:3Sasak}

We cover some of the basics of 3-Sasakian manifolds and 3-Sasakian reduction.  See~\cite{BG2} for more details.
\begin{defn}\label{defn:3Sasak}
Let $(\mathcal{S},g)$ be a Riemannian manifold of dimension $n=4m+3$.  Then $\mathcal{S}$ is 3-Sasakian
if it admits three Killing vector fields $\{\xi^1,\xi^2,\xi^3\}$ each satisfying definition (\ref{defn:Sasak})
such that $g(\xi^i,\xi^k)=\delta_{ij}$ and $[\xi^i,\xi^j]=2\epsilon_{ijk}\xi^k$.
\end{defn}
We have a triple of Sasakian structures on $\mathcal{S}$.  For $i=1,2,3$ we have $\eta^i(X)=g(\xi^i,X)$ and
$\Phi^i(X)=\nabla_X \xi^i$.  We say that $\{g,\Phi^i,\xi^i,\eta^i:i=1,2,3\}$ defines a \emph{3-Sasakian
structure} on $\mathcal{S}$.

\begin{prop}\label{prop:quat-ident}
The tensors $\Phi^i,i=1,2,3$ satisfy the following identities.
\begin{thmlist}
\item $\Phi^i(\xi^j)=-\epsilon_{ijk}\xi^k,$
\item $\Phi^i\circ\Phi^j =-\epsilon_{ijk}\Phi^k +\xi^i\otimes\eta^j -\delta_{ij}Id$
\end{thmlist}
\end{prop}

Notice that if $\alpha=(a_1,a_2,a_3)\in S^2\subset\R^3$ then
$\xi(\alpha)=a_1\xi^1+a_2\xi^2+a_3\xi^3$ is a Sasakian structure.  Thus a
3-Sasakian manifold come equipped with an $S^2$ of Sasakian structures.

As in definition~\ref{defn:Sasak}, 3-Sasakian manifolds can be characterized by the holonomy of the cone $C(\mathcal{S})$.
\begin{prop}\label{prop:3Sasak-cone}
Let $(\mathcal{S},g)$ be a Riemannian manifold of dimension $n=4m+3$.  Then $(\mathcal{S},g)$ is
3-Sasakian if, and only if, the holonomy of the metric cone $(C(\mathcal{S}),\ol{g})$ is a subgroup of
$Sp(m+1)$.  In other words, $(C(\mathcal{S}),\ol{g})$ is hyperk\"{a}hler.
\end{prop}
\begin{proof}
Define almost complex structures $I_i,i=1,2,3$ by
\[I_i X =-\Phi^i(X)+\eta^i(X)r\partial_r, \text{  and }I_i r\partial_r=\xi^i. \]
It is straight forward to verify that they satisfy
$I_i\circ I_j = \epsilon_{ijk}I_k -\delta_{ij}Id$.  And from the
integrability condition on each $\Phi^i$ in definition (\ref{defn:Sasak}) each
$I_i,i=1,2,3$ is parallel.
\end{proof}

Since a hyperk\"{a}hler manifold is Ricci flat, we have the following.
\begin{cor}
A 3-Sasakian manifold $(\mathcal{S},g)$ of dimension $n=4m+3$ is Einstein with positive scalar
curvature $s=2(2m+1)(4m+3)$.  Furthermore, if $(\mathcal{S},g)$ is complete, then it is compact
with finite fundamental group.
\end{cor}

The structure group of a 3-Sasakian manifold reduces to $Sp(m)\times\mathbb{I}_3$ where $\mathbb{I}_3$
is the $3\times 3$ identity matrix.  Thus we have
\begin{cor}
A 3-Sasakian manifold $(M,g)$ is spin.
\end{cor}

Suppose $(\mathcal{S},g)$ is compact.  This will be the case in all examples considered here.
Then the vector fields $\{\xi^1,\xi^2,\xi^3\}$ are complete and define a locally free action
of $Sp(1)$ on $(\mathcal{S},g)$.  This defines a foliation $\mathscr{F}_3$,
the \emph{3-Sasakian foliation}.  The generic leaf is either $SO(3)$ or $Sp(1)$, and
all the leaves are compact.  So $\mathscr{F}_3$ is quasi-regular, and the space
of leaves is a compact orbifold, denoted $\mathcal{M}$.  The projection
$\varpi:\mathcal{S}\rightarrow\mathcal{M}$ exhibits $\mathcal{S}$ as an $SO(3)$ or $Sp(1)$
V-bundle over $\mathcal{M}$.  The leaves of $\mathscr{F}_3$ are constant curvature
3-Sasakian 3-manifolds which must be homogeneous spherical space forms.  Thus a leaf is
$\Gamma\backslash S^3$ with $\Gamma\subset Sp(1)$.

For $\beta\in S^2$ we also have the characteristic vector field $\xi_\beta$ with the associated 1-dimensional
foliation $\mathscr{F}_\beta \subset\mathscr{F}_3$.  In this case $\mathscr{F}_\beta$ is automatically
quasi-regular. Denote the leaf space of $\mathscr{F}_\beta$ as $\mathcal{Z}_\beta$ or just $\mathcal{Z}$. Then
the natural projection $\pi:\mathcal{S}\rightarrow\mathcal{Z}$ is an $S^1$ V-bundle.  And
$\mathcal{Z}$ is a $(2m+1)$-dimensional projective, normal algebraic variety with orbifold singularities
and a K\"{a}hler form $\omega\in c_1^{orb}(\mathcal{Z})$, i.e. is Fano.

Fix a Sasakian structure $\{\Phi^1,\xi^1,\eta^1\}$ on $\mathcal{S}$.  The horizontal subbundle
$\mathcal{H}=\ker\eta^1$ to the foliation $\mathcal{F}$ of $\xi^1$ with the almost complex structure
$I=-\Phi^1|_{\mathcal{H}}$ define a CR structure on $\mathcal{S}$.  The form
$\eta=\eta^2+i\eta^3$ is of type $(1,0)$ with respect to $I$.
And $d\eta|_{\mathcal{H}\cap\ker(\eta)}\in\Omega^{2,0}(\mathcal{H}\cap\ker(\eta))$ is nondegenerate as a complex 2-form on
$\mathcal{H}\cap\ker(\eta)$.  Consider the complex 1-dimensional subspace $P\subset\Lambda^{1,0}\mathcal{H}$
spanned by $\eta$.
Letting $\exp(it\xi^1)$ denote an element of the circle subgroup $U(1)\subset Sp(1)$ generated by
$\xi^1$ one see that $\exp(it\xi^1)$ acts on $P$ with character $e^{-2it}$.  Then
$\mathbf{L}\cong\mathcal{S}\times_{U(1)}P$ defines a holomorphic line V-bundle over
$\mathcal{Z}$.  And we have a holomorphic section $\theta$ of $\Lambda^{1,0}(\mathcal{Z})\otimes\mathbf{L}$
such that
\[ \theta(X)=\eta(\tilde{X}), \]
where $\tilde{X}$ is the horizontal lift of a vector field $X$ on $\mathcal{Z}$.
Let $D=\ker(\theta)$ be the complex distribution defined by $\theta$.  Then
$d\theta|_D \in\Gamma(\Lambda^2 D\otimes\mathbf{L})$ is nondegenerate.  Thus $D=\ker(\theta)$ is
complex contact structure on $\mathcal{Z}$, that is, a maximally non-integrable holomorphic subbundle
of $T^{1,0}\mathcal{Z}$.  Also, $\theta\wedge (d\theta)^m$ is a nowhere zero section of
$\mathbf{K}_{\mathcal{Z}}\otimes\mathbf{L}^{m+1}$.  Thus
$\mathbf{L}\cong\mathbf{K}_{\mathcal{Z}}^{-\frac{1}{m+1}}$ as holomorphic line V-bundles.
We have the following for 3-Sasakian manifolds.

\begin{thm}\label{thm:leaf-3Sasak}
Let $(\mathcal{S},g)$ be a compact 3-Sasakian manifold of dimension $n=4m+3$, and let $\mathcal{Z}_\beta$ be
the leaf space of the foliation $\mathcal{F}_\beta$ for $\beta\in S^2$.  Then $\mathcal{Z}_\beta$ is
a compact $\Q$-factorial contact Fano variety with a K\"{a}hler-Einstein metric $h$ with scalar curvature
$s=8(2m+1)(m+1)$.  The projection $\pi:\mathcal{S}\rightarrow\mathcal{Z}$ is an orbifold Riemannian
submersion with respect to the metrics $g$ on $\mathcal{S}$ and $h$ on $\mathcal{Z}$.
\end{thm}

The space $\mathcal{Z}=\mathcal{Z}_\beta$ is, up to isomorphism of all structures, independent of
$\beta\in S^2$.  We call $\mathcal{Z}$ the \emph{twistor space} of $\mathcal{S}$.
Consider again the natural projection $\varpi:\mathcal{S}\rightarrow\mathcal{M}$ coming from
the foliation $\mathscr{F}_3$.  This factors into $\pi:\mathcal{S}\rightarrow\mathcal{Z}$
and $\rho:\mathcal{Z}\rightarrow\mathcal{M}$.  The generic fibers of $\rho$ is a $\cps^1$ and there are
possible singular fibers $\Gamma\backslash\cps^1$ which are simply connected and for which $\Gamma\subset U(1)$
is a finite group.  And restricting to a fiber $\mathbf{L}|_{\cps^1}=\mathcal{O}(2)$, which is a
V-bundle on singular fibers.
Consider $g=\exp(\frac{\pi}{2}\xi^2)\in Sp(1)$ which gives an isometry of $\mathcal{S}$
$\varsigma_g :\mathcal{S}\rightarrow\mathcal{S}$ for which $\varsigma_g(\xi^1)=-\xi^1$.  And $\varsigma_g$
descends to an anti-holomorphic isometry $\sigma:\mathcal{Z}\rightarrow\mathcal{Z}$ preserving the fibers.

We now consider the orbifold $\mathcal{M}$ more closely.  Let $(\mathcal{M},g)$ be any $4m$ dimensional
Riemannian orbifold.  An \emph{almost quaternionic} structure on $\mathcal{M}$ is
a rank 3 V-subbundle $\mathcal{Q}\subset End(T\mathcal{M})$ which is locally spanned by almost complex
structures $\{J_i\}_{i=1,2,3}$ satisfying the quaternionic identities $J_i^2=-Id$ and
$J_1 J_2 =-J_2 J_1 =J_3$.  We say that $\mathcal{Q}$ is compatible with $g$ if $J_i^*g=g$ for $i=1,2,3$.
Equivalently, each $J_i,i=1,2,3$ is skew symmetric.
\begin{defn}\label{defn:quatern-kahler}
A Riemannian orbifold $(\mathcal{M},g)$ of dimension $4m, m>1$ is \emph{quaternionic K\"{a}hler} if there is an almost
quaternionic structure $\mathcal{Q}$ compatible with $g$ which is preserved by the Levi-Civita
connection.
\end{defn}
This definition is equivalent to the holonomy of $(\mathcal{M},g)$ being contained in
\linebreak $Sp(1)Sp(m)$.
For orbifolds this is the holonomy on $\mathcal{M}\setminus S_{\mathcal{M}}$ where $S_{\mathcal{M}}$ is
the singular locus of $\mathcal{M}$.  For more on quaternionic K\"{a}hler manifolds see~\cite{S1}.
Notice that this definition always holds on an oriented Riemannian 4-manifold ($m=1$).
This case requires a different definition.  Consider the \emph{curvature operator}
$\mathcal{R}:\Lambda^2\rightarrow\Lambda^2$ of an oriented Riemannian 4-manifold.
With respect to the decomposition $\Lambda^2=\Lambda^2_+ \oplus\Lambda^2_-$, we have
\begin{equation}\label{mat:curv-dec}
\mathcal{R}=
\left\lgroup
\begin{matrix}
W_{+}+\frac{s}{12} & \overset{\circ}{r} \\
\overset{\circ}{r} & W_{-}+\frac{s}{12}
\end{matrix}\right\rgroup,
\end{equation}
where $W_+$ and $W_-$ are the selfdual and anti-self-dual pieces of the Weyl curvature and
$\overset{\circ}r=\Ric -\frac{s}{4}g$ is the trace-free Ricci curvature.
An oriented 4 dimensional Riemannian orbifold $(\mathcal{M},g)$ is quaternionic K\"{a}hler if it is Einstein
and anti-self-dual, meaning that $\overset{\circ}{r}=0$ and $W_+ =0$.
\begin{thm}\label{thm:QK-3Sasak}
Let $(\mathcal{S},g)$ be a compact 3-Sasakian manifold of dimension $n=4m+3$.  Then there is a natural
quaternionic K\"{a}hler structure on the leaf space of $\mathcal{F}_3$, $(\mathcal{M},\check{g})$, such that
the V-bundle map $\varpi:\mathcal{S}\rightarrow\mathcal{M}$ is a Riemannian submersion.
Furthermore, $(\mathcal{M},\check{g})$ is Einstein with scalar curvature $16m(m+2)$.
\end{thm}

\subsection{3-Sasakian reduction}

We now summarize 3-Sasakian reduction and its application to producing infinitely many 3-Sasakian 7-manifolds.
In particular, we are interested in toric 3-Sasakian 7-manifolds which
have a $T^2$ action preserving the 3-Sasakian structure.  Up to coverings they are all
obtainable by taking 3-Sasakian quotients of $S^{4n-1}$ by a torus $T^{k}$, $k=n-2$.
See~\cite{BGMR,BG2} for more details.

Let $(\mathcal{S},g)$ be a 3-Sasakian manifold.  And let $I(\mathcal{S},g)$ be the subgroup in the isometry
group $\Isom(\mathcal{S},g)$ of 3-Sasakian automorphisms.
\begin{defn}
Let $(\mathcal{S},g)$ be a 3-Sasakian 7-manifold.  Then $(\mathcal{S},g)$ is \emph{toric} if there is a real
2-torus $T^2\subseteq I(\mathcal{S},g)$.
\end{defn}

Let $G\subseteq I(\mathcal{S},g)$ be compact with Lie algebra $\mathfrak{g}$.
One can define the \emph{3-Sasakian moment map}
\begin{equation}
\mu_{\mathcal{S}}:\mathcal{S}\rightarrow\mathfrak{g}^*\otimes\R^3
\end{equation}
by
\begin{equation}
\langle\mu^a_{\mathcal{S}}, X\rangle = \frac{1}{2}\eta^a(\tilde{X}),\quad  a=1,2,3\text{  for }X\in\mathfrak{g},
\end{equation}
where $\tilde{X}$ be the vector field on $\mathcal{S}$ induced by $X\in\mathfrak{g}$.

\begin{prop}\label{prop:3-Sasak-red}
Let $(\mathcal{S},g)$ be a 3-Sasakian manifold and $G\subset I(\mathcal{S},g)$ a connected compact
subgroup.  Assume that $G$ acts freely (locally freely) on $\mu_\mathcal{S}^{-1}(0)$.
Then $\mathcal{S}\sslash G = \mu_\mathcal{S}^{-1}(0)/G$ has the structure of a 3-Sasakian manifold (orbifold).
Let $\iota:\mu_{\mathcal{S}}^{-1}(0)\rightarrow\mathcal{S}$ and
$\pi:\mu_{\mathcal{S}}^{-1}(0)\rightarrow\mu_{\mathcal{S}}^{-1}(0)/G$ be the corresponding
embedding and submersion.  Then the metric $\check{g}$ and 3-Sasakian vector fields are defined by
$\pi^*\check{g}=\iota^*g$ and $\pi_*\xi^i|_{\mu_\mathcal{S}^{-1}(0)}=\check{\xi}^i$.
\end{prop}

Consider the unit sphere $S^{4n-1}\subset\Ha^n$ with the metric $g$ obtained by restricting the flat metric on
$\Ha^n$.  Give $S^{4n-1}$ the standard 3-Sasakian structure induced by the right action of $Sp(1)$.
Then $I(S^{4n-1},g)=Sp(n)$ acting by the standard linear representation on the left.
We have the maximal torus $T^n\subset Sp(n)$ and every representation of a subtorus
$T^k$ is conjugate to an inclusion $\iota_\Omega :T^k\rightarrow T^n$ which is
represented by a matrix
\begin{equation}\label{eq:torus-inc}
\iota_\Omega(\tau_1,\ldots ,\tau_k)=
\left\lgroup
\begin{matrix}
\prod_{i=1}^{k}\tau_1^{a_1^i} & \cdots & 0 \\
\vdots    & \ddots & \vdots \\
0         & \cdots & \prod_{i=1}^k\tau_n^{a_n^i} \\
\end{matrix}\right\rgroup,
\end{equation}
where $(\tau_1,\ldots, \tau_k)\in T^k$.  Every such representation is defined by the
$k\times n$ integral \emph{ weight matrix}
\begin{equation}
\Omega=\left\lgroup
\begin{matrix}
a_1^1 & \cdots & a_k^1 & \cdots & a_n^1 \\
a_1^2 & \cdots & a_k^2 & \cdots & a_n^2 \\
\vdots & \ddots & \vdots & \ddots & \vdots \\
a_1^k & \cdots & a_k^n & \cdots & a_n^k \\
\end{matrix}\right\rgroup
\end{equation}

Let $\{e_i\},i=1,\ldots,k$ be a basis for the dual of the Lie algebra of $T^k$, $\mathfrak{t}_k^*\cong\R^k$.
Then the moment map $\mu_\Omega:S^{4n-1}\rightarrow\mathfrak{t}_k^*\otimes\R^3$ can be written
as $\mu_\Omega =\sum_j\mu_\Omega^j e_j$, where in terms of complex coordinates $u_l=z_l+w_l j$ on $\Ha^n$
\begin{equation}\label{eq:3-Sasak-mom}
\mu_\Omega^j(\mathbf{z},\mathbf{w})=i\sum_l a_l^j(|z_l|^2 - |w_l|^2)+2k\sum_l a_l^j\ol{w}_lz_l.
\end{equation}

Denote by $\Delta_{\alpha_1,\ldots,\alpha_k}$ the $\binom{n}{k}$ $k\times k$ minor determinants of $\Omega$.
\begin{defn}\label{defn:adm}
Let $\Omega\in\mathscr{M}_{k,n}(\Z)$ be a weight matrix.
\begin{thmlist}
\item $\Omega$ is \emph{nondegenerate} if $\Delta_{\alpha_1,\ldots,\alpha_k}\neq 0$,
for all $1\leq\alpha_1 < \cdots < \alpha_k \leq n$.\\

Let $\Omega$ be nondegenerate, and let $d$ be the $\gcd$ of all the $\Delta_{\alpha_1,\ldots,\alpha_k}$,
the kth determinantal divisor.  Then $\Omega$ is \emph{admissible} if
\item $\gcd(\Delta_{\alpha_2,\ldots,\alpha_{k+1}},\ldots,\Delta_{\alpha_1,\ldots,\hat{\alpha}_t,\ldots,\alpha_{k+1}},
\ldots,\Delta_{\alpha_1,\ldots,\alpha_k})=d$ for all length $k+1$ sequences
$1\leq\alpha_1 <\cdots <\alpha_t <\cdots <\alpha_{k+1}\leq n+1$.
\end{thmlist}
\end{defn}

The quotient obtained in Proposition~\ref{prop:3-Sasak-red}
$\mathcal{S}_\Omega =S^{4n-1}\sslash T^k(\Omega)$ will depend on $\Omega$ only up
to a certain equivalence.  Choosing a different basis of $\mathfrak{t}_k$ results in an
action on $\Omega$ by an element in $Gl(k,\Z)$.  We also have the normalizer of $T^n$ in
$Sp(n)$, the Weyl group $\mathscr{W}(Sp(n))=\Sigma_n \times\Z_2^n$ where $\Sigma_n$ is the permutation group.
$\mathscr{W}(Sp(n))$ acts on $S^{4n-1}$ preserving the 3-Sasakian structure, and it acts on
weight matrices by permutations and sign changes of columns.
The group $Gl(k,\Z)\times\mathscr{W}(Sp(n))$ acts on $\mathscr{M}_{k,n}(\Z)$.
We say $\Omega$ is \emph{reduced} if $d=1$ in definition~\ref{defn:adm}.  It is a result in~\cite{BGMR} that we may
assume that a nondegenerate weight matrix $\Omega$ is reduced, as this is the case when (\ref{eq:torus-inc})  is an inclusion.
\begin{thm}[\cite{BG2,BGMR}]
Let $\Omega\in\mathscr{M}_{k,n}(\Z)$ be reduced.
\begin{thmlist}
\item  If $\Omega$ is nondegenerate, then $\mathcal{S}_\Omega$ is an orbifold.\\
\noindent
\item  Supposing $\Omega$ is nondegenerate, $\mathcal{S}_\Omega$ is smooth if and only if
$\Omega$ is admissible.
\end{thmlist}
\end{thm}
Notice that the automorphism group of $\mathcal{S}_\Omega$ contains
$T^{n-k}\cong T^n/\iota_\Omega(T^k)$.

We now restrict to the case of 7-dimensional toric quotients, so $n=k+2$.  We may take
matrices of the form
\begin{equation}\label{eq:3-Sasak-mat}
\Omega =\left\lgroup
\begin{matrix}
1 & 0 & \cdots & 0 & a_1 & b_1 \\
0 & 1 & \cdots & 0 & a_2 & b_2 \\
\vdots & \vdots & \ddots & \vdots & \vdots & \vdots \\
0 & 0 & \cdots & 1 & a_k & b_k
\end{matrix}
\right\rgroup.
\end{equation}

\begin{prop}[\cite{BGMR}]\label{prop:3-Sasak-mat}
Let $\Omega\in\mathscr{M}_{k,k+2}(\Z)$ be as above.  Then $\Omega$ is admissible if and only if
$a_i,b_j,i,j=1,\ldots, k$ are all nonzero, $\gcd(a_i,b_i)=1$ for $i=1,\ldots,k$, and we do not have
$a_i= a_j$ and $b_i=b_j$, or $a_i=-a_j$ and $b_i=-b_j$ for some $i\neq j$.
\end{prop}

This shows that for $n=k+2$ there are infinitely many reduced admissible
weight matrices.  One can, for example, choose $a_i,b_j,i,j=1,\ldots k$ be all pairwise relatively prime.
We will make use of the cohomology computation of R. Hepworth~\cite{He} to show that we have infinitely many
smooth 3-Sasakian 7-manifolds of each second Betti number $b_2\geq 1$.
Let $\Delta_{p,q}$ denote the $k\times k$ minor determinant of $\Omega$ obtained by deleting the $p^{th}$ and $q^{th}$
columns.

\begin{thm}[\cite{He}\cite{BGM2,BGMR}]\label{thm:3-Sasak-coh}
Let $\Omega\in\mathscr{M}_{k,k+2}(\Z)$ be a reduced admissible weight matrix.  Then
$\pi_1(\mathcal{S}_\Omega)=e$.  And the cohomology of $\mathcal{S}_\Omega$ is

\renewcommand{\arraystretch}{1.5}
\begin{centering}
\begin{tabular}{|l|cccccccc|}\hline
$p$ & $0$ & $1$ & $2$ & $3$ & $4$ & $5$ & $6$ & $7$\\ \hline
$H^p$ & $\Z$ & $0$ & $\Z^k$ & $0$ & $G_\Omega$ & $\Z^k$ & $0$ & $\Z$ \\\hline
\end{tabular},\\
\end{centering}
where $G_\Omega$ is a torsion group of order
\[ \sum |\Delta_{s_1, t_1}|\cdots |\Delta_{s_{k+1},t_{k+1}}| \]
with the summand with index $s_1,t_1,\ldots, s_{k+1},t_{k+1}$ included if and only if
the graph on the vertices $\{1,\ldots, k+2\}$ with edges $\{s_i,t_i\}$ is a tree.
\end{thm}

If we consider weight matrices as in Proposition~\ref{prop:3-Sasak-mat}, then the order
of $G_\Omega$ is greater than $|a_1\cdots a_k|+|b_1\cdots b_k|$.
We have the following.
\begin{cor}[\cite{He}\cite{BGMR}]\label{cor:3-Sasak-coh}
There are smooth toric 3-Sasakian 7-manifolds with second Betti number $b_2=k$ for all $k\geq 0$.
Furthermore, there are infinitely many possible homotopy types of examples $\mathcal{S}_{\Omega}$ for each
$k>0$.
\end{cor}

\subsection{Anti-self-dual Einstein spaces}\label{sec:asd-Einst}
\sloppy

We consider the anti-self-dual Einstein orbifolds $\mathcal{M}=\mathcal{M}_\Omega$ associated to the toric 3-Sasakian 7-manifolds
$\mathcal{S}$ in greater detail.  Since $\mathcal{M}$ is a 4 dimensional orbifold with an effective action of
$T^2$ the techniques of~\cite{HS2} show that $\mathcal{M}$ is characterized by the polygon $\mathcal{Q}_\Omega =\mathcal{M}/T^2$
with $k+2$ edges, $b_2(\mathcal{M})=k$, labeled in cyclic order with
$(m_0,n_0),(m_1,n_1),\ldots, (m_{k+2},n_{k+2})$ in $\Z^2$, $(m_0,n_0)=-(m_{k+2},n_{k+2})$, denoting the isotropy subgroups.
For the quotients $\mathcal{M}_\Omega$ one can show the following (cf~\cite{vC2}):
\renewcommand{\theenumi}{\alph{enumi}}
\begin{enumerate}
\renewcommand{\labelenumi}{\theenumi.}
\item  The sequence $m_i, \ i=0,\ldots k+2$ is strictly increasing.

\item  The sequence $(n_i - n_{i-1})/(m_i -m_{i-1}), \ i=1,\ldots k+2$ is strictly increasing.
\end{enumerate}
We will make use of the following result of D. Calderbank and M. Singer~\cite{CaSin} which classifies those
compact orbifolds which admit toric anti-self-dual Einstein metrics.
The case for which the associated 3-Sasakian space is smooth is originally due to R. Bielawski~\cite{Bie}.
\begin{thm}\label{thm:asd-class}
Let $\mathcal{M}$ be a compact toric 4-orbifold with $\pi^{orb}_1(\mathcal{M})=e$ and $k=b_2(\mathcal{M})$.
Then the following are equivalent.
\begin{thmlist}
\item  One can arrange that the isotropy data of $\mathcal{M}$ satisfy a.\ and b.\ above by
cyclic permutations, changing signs, and acting by $Gl(2,\Z)$.
\item  $\mathcal{M}$ admits a toric anti-self-dual Einstein metric unique up to homothety and
equivariant diffeomorphism.  Furthermore, $(\mathcal{M},g)$ is isometric to the quaternionic K\"{a}hler
reduction of $\qps^{k+1}$ by a torus $T^k\subset Sp(k+2)$.
\end{thmlist}
\end{thm}
It is well known that the only possible smooth compact anti-self-dual Einstein spaces with positive scalar curvature
are $S^4$ and $\overline{\cps}^2$, which are both toric.

Suppose $\mathcal{M}$ has isotropy data $v_0,v_1,\ldots,v_{k+2}$.  Then it is immediate that
$v_0,v_1,\ldots,v_{k+2},-v_1,-v_2,\ldots,-v_{k+1}$ are the vertices of a convex polygon
in $N_\R =\R^2$, which defines an augmented fan $\Delta^*$ defining a toric Fano surface $X$.
A bit more thought gives the following.
\begin{thm}
There is a one to one correspondence between compact toric anti-self-dual Einstein orbifolds $\mathcal{M}$ with
$\pi^{orb}_1(\mathcal{M})=e$ and
special symmetric toric Fano orbifold surfaces $X$ with $\pi^{orb}_1(X)=e$.
And $X$ has a K\"{a}hler-Einstein metric of positive scalar curvature.
Under the correspondence if $b_2(\mathcal{M})=k$, then $b_2(X)=2k+2$.
\end{thm}
This will be reproved in section~\ref{sec:twist-div} by exhibiting $X$ as a divisor in the twistor space.

\begin{xpl}
Consider the admissible weight matrix
\[\Omega=\left\lgroup
\begin{matrix}
1 & 0 & 1 & 1 \\
0 & 1 & 1 & 2
\end{matrix}\right\rgroup.\]
Then the 3-Sasakian space $\mathcal{S}_\Omega$ is smooth and
$b_2(\mathcal{S}_\Omega)=b_2(\mathcal{M}_\Omega)=2$.
And the anti-self-dual orbifold $\mathcal{M}_\Omega$ has isotropy data
\[v_0 =(-7,-2),(-5,-2),(-1,-1),(5,1),(7,2)=v_4.\]
The singular set of $\mathcal{M}$ consists of two points with stablizer group $\Z_3$ and two with
$\Z_4$.  The associated toric K\"{a}hler-Einstein surface is that in figure \ref{fig:typ-xpl}.
\end{xpl}

The generic fiber of $\varpi:\mathcal{S}\rightarrow\mathcal{M}$ is either $SO(3)$ of $Sp(1)$.
In general the existence of a lifting to an $Sp(1)$ V-bundle is obstructed by the
\emph{Marchiafava-Romani class} which when $\mathcal{M}$ is 4-dimensional is identical to
$w_2(\mathcal{M})\in H^2_{orb}(\mathcal{M},\Z_2)$.  In other words, the contact line bundle $\mathbf{L}$ has
a square root $\mathbf{L}^{\frac{1}{2}}$ if, and only if, $w_2(\mathcal{M})=0$.
\begin{prop}
Let $X$ be the symmetric toric Fano surface associated to the anti-self-dual Einstein orbifold $\mathcal{M}$.
Then $\ind(X)=2$ if and only if $w_2(\mathcal{M})=0$.
In other words, $\mathbf{K}^{-1}_X$ has a square root if and only if the contact line bundle $\mathbf{L}$ does.
\end{prop}
See~\cite{vC2} Proposition 5.22 for a proof.

\subsection{Twistor space and divisors}\label{sec:twist-div}

We will consider the twistor space $\mathcal{Z}$ introduced in Theorem~\ref{thm:leaf-3Sasak}
more closely for the case when $\mathcal{M}$ is an anti-self-dual Einstein orbifold.
For now suppose $(\mathcal{M},[g])$ is an anti-self-dual, i.e. $W_+ \equiv 0$,
conformal orbifold.  There is a twistor space of $(\mathcal{M},[g])$ which is originally due to R. Penrose
~\cite{Pen}.  See~\cite{AHS} for positive definite case.

The \emph{twistor space} of $(\mathcal{M},[g])$ is a complex three dimensional orbifold $\mathcal{Z}$
with the following properties:
\renewcommand{\theenumi}{\alph{enumi}}
\begin{enumerate}
\renewcommand{\labelenumi}{\theenumi.}
\item  There is a $V$-bundle fibration $\varpi:\mathcal{Z}\rightarrow\mathcal{M}$.
\item  The general fiber of $P_x =\varpi^{-1}(x), x\in\mathcal{Z}$ is a projective line
$\cps^1$ with normal bundle $N\cong\mathcal{O}(1)\oplus\mathcal{O}(1)$, which holds over singular fibers
with $N$ a V-bundle.
\item  There exists an anti-holomorphic involution $\sigma$ of $\mathcal{Z}$ leaving the fibers
$P_x$ invariant.
\end{enumerate}

Let $T$ be an oriented real 4-dimensional vector space with inner product $g$.  Let
$C(T)$ be set of orthogonal complex structures inducing the orientation, i.e.
if $r,s\in T$ is a complex basis then $r,Jr,s,Js$ defines the orientation.  One has
$C(T)= S^2\subset\Lambda^2_+(T)$, where $S^2$ is the sphere of radius $\sqrt{2}$.
Now take $T$ to be $\Ha$.
Recall that $Sp(1)$ is the group of unit quaternions.  Let
\begin{equation}
Sp(1)_+ \times Sp(1)_-
\end{equation}
act on $\Ha$ by
\begin{equation}
w\rightarrow gw{g^\prime}^{-1},\text{  for  }w\in\Ha\text{ and } (g,g^\prime)\in Sp(1)_+ \times Sp(1)_-.
\end{equation}
Then we have
\begin{equation}
Sp(1)_+ \times_{\Z_2} Sp(1)_- \cong SO(4),
\end{equation}
where $\Z_2$ is generated by $(-1,-1)$.
Let
\begin{equation}
\begin{split}
C & =\{ai+bj+ck: a^2+b^2+c^2=1, a,b,c\in\R\}\\
  & =\{g\in Sp(1)_+: g^2=-1\}\cong S^2.
\end{split}
\end{equation}
Then $g\in C$ defines an orthogonal complex structure by
\[ w\rightarrow gw,\text{  for  } w\in\Ha,\]
giving an identification $C=C(\Ha)$.
Let $V_+ =\Ha$ considered as a representation of $Sp(1)_+$ and a right $\C$-vector space.
Define $\pi:V_+\setminus\{0\}\rightarrow C$ by $\pi(h)=-hih^{-1}$.  Then the fiber of
$\pi$ over $hih^{-1}$ is $h\C$.  Then $\pi$ is equivariant if $Sp(1)_+$ acts on $C$ by
$q\rightarrow gqg^{-1}, g\in Sp(1)_+$.  We have a the identification
\begin{equation}
C=V_+\setminus\{0\}/{\C^*} =\mathbb{P}(V_+).
\end{equation}

Fix a Riemannian metric $g$ in $[g]$.  Let $\phi:\tilde{U}\rightarrow U\subset\mathcal{M}$
be a local uniformizing chart with group $\Gamma$.
Let $F_{\tilde{U}}$ be the bundle of orthonormal frames on $\tilde{U}$.
Then
\begin{equation}\label{eq:twistor-chart}
F_{\tilde{U}}\times_{SO(4)}\mathbb{P}(V_+) =F_{\tilde{U}}\times_{SO(4)} C
\end{equation}
defines a local uniformizing chart for $\mathcal{Z}$ mapping to
\[F_{\tilde{U}}\times_{SO(4)}\mathbb{P}(V_+)/\Gamma =F_{\tilde{U}}/\Gamma \times_{SO(4)}\mathbb{P}(V_+).\]
Right multiplication by $j$ on $V_+ =\Ha$
defines the anti-holomorphic involution $\sigma$ which is fixed point free on (\ref{eq:twistor-chart}).
We will denote a neighborhood as in (\ref{eq:twistor-chart}) by $\tilde{U}_{\mathcal{Z}}$.

An almost complex structure is defined as follows.  At a point $z\in\tilde{U}_{\mathcal{Z}}$ the Levi-Civita
connection defines a horizontal subspace $H_z$ of the real tangent space $T_z$ and we have a
splitting
\begin{equation}\label{eq:twistor-comp}
T_z =H_z\oplus T_z P_x = T_x\oplus T_z P_x,
\end{equation}
where $\varpi(z)=x$ and $T_x$ is the real tangent space of $\tilde{U}$.
Let $J_z$ be the complex structure on $T_x$ given by $z\in P_x=C(T_x)$, and let
$J_z^\prime$ be complex structure on $T_x\oplus T_z P_x$ arising from the natural complex structure
on $P_x$.  Then the almost complex structure on $T_z$ is the direct sum of $J_z$ and $J_z^\prime$.
This defines a natural almost complex structure on $Z_{\tilde{U}}$ which is invariant under $\Gamma$.
We get an almost complex structure on $\mathcal{Z}$ which is integrable precisely when $W_+\equiv 0$.

Assume that $\mathcal{M}$ anti-self-dual Einstein with non-zero scalar curvature.  Then
$\mathcal{Z}$ has a complex contact structure $D\subset T^{1,0}\mathcal{Z}$ with
holomorphic contact form $\theta\in\Gamma(\Lambda^{1,0}\mathcal{Z}\otimes\mathbf{L})$ where
$\mathbf{L}=T^{1,0}\mathcal{Z}/D$.

The group of isometries $\Isom(\mathcal{M})$ lifts to an action on $\mathcal{Z}$ by real holomorphic transformations.
Real means commuting with $\sigma$.  This extends to a holomorphic action of the complexification
$\Isom(\mathcal{M})_{\C}$.
For $X\in\mathfrak{Isom}(\mathcal{M})\otimes\C$, the Lie algebra of $\Isom(\mathcal{M})_{\C}$,
we will also denote by $X$
the holomorphic vector field induced on $\mathcal{Z}$.
Then $\theta(X)\in H^0(\mathcal{Z},\mathcal{O}(\mathbf{L}))$.
By a well known twistor correspondence the map $X\rightarrow\theta(X)$ defines an isomorphism
\begin{equation}\label{eq:twistor-trans}
\mathfrak{Isom}(\mathcal{M})\otimes\C\cong H^0(\mathcal{Z},\mathcal{O}(\mathbf{L})),
\end{equation}
which maps real vector fields to real sections of $\mathbf{L}$.

Suppose for now on that $\mathcal{M}$ is a toric anti-self-dual Einstein orbifold with twistor space $\mathcal{Z}$.
We will assume that $\pi_1^{orb}(\mathcal{M})=e$ which can always be arranged by taking the orbifold cover.
Then as above $T^2$ acts on $\mathcal{Z}$ by holomorphic transformations.  And the action extends to
$T^2_\C =\C^*\times\C^*$, which in this case is an algebraic action.
Let $\mathfrak{t}$ be the Lie algebra of $T^2$ with $\mathfrak{t}_\C$ the Lie algebra of $T^2_\C$.  Then we have
from (\ref{eq:twistor-trans}) the pencil
\begin{equation}
E=\mathbb{P}(\mathfrak{t}_\C)\subseteq |\mathbf{L}|,
\end{equation}
where for $t\in E$ we denote $X_t=(\theta(t))$ the divisor of the section $\theta(t)\in H^0(\mathcal{Z},\mathcal{O}(\mathbf{L}))$.
Note that $E$ has an equator of real divisors.  Also, since $T^2_\C$ is abelian, every $X_t, t\in E$ is
$T^2_\C$ invariant.

Our goal is to determine the structure of the divisors in the pencil $E$.
As before we will consider the one parameter groups $\rho_i \in N=\Z\times\Z$, where $N$ is the lattice
of one parameter $\C^*$-subgroups of $T^2_\C$.
Also, we will identify the Lie algebra $\mathfrak{t}$ of $T^2$ with $N\otimes\R$ and
the Lie algebra $\mathfrak{t}_\C$ of $T^2_\C$ with $N\otimes\C$.
Since $\mathbf{L}|_{P_x} =\mathcal{O}(2)$ a divisor $X_t\in E$ intersects a generic twistor line
$P_x$ at two points.

Recall that the set of non-trivial stablizers of the $T^2$-action on $\mathcal{M}$ is
$B=\cup_{i=1}^{k+2} B_i$ where $B_i$ is topologically a 2-sphere.  Denote by $x_i =B_i\cap B_{i+1}$
the $k+2$ fixed points of the action.
We will denote $P_i :=P_{x_i}, i=1,\ldots, k+2$.
One can show there exist two irreducible rational curves $C_i^{\pm}, i=1,\ldots, k+2$ mapped diffeomorphicly to
$B_i$ by $\varpi$.  Furthermore, $\sigma(C_i^{\pm})=C_i^{\mp}$.
The singular set for the $T^2$-action on $\mathcal{Z}$ is the union of rational curves
\begin{equation}\label{eq:twist-sing}
\Sigma=\Bigl(\cup_{i=1}^{k+2}P_i \Bigr)\bigcup\Bigl(\cup_{i+1}^{k+2}C_i^+\cup C_i^-\Bigr).
\end{equation}

With a closer analysis of the action of $T_\C^2$ on $\mathcal{Z}$ one can prove the following.
See~\cite{vC2} for the proof.
The term \emph{suborbifold} denotes a subvariety which is a submanifold in every local uniformizing neighborhood.

\begin{thm}\label{thm:twistor-divisor}
Let $\mathcal{M}$ be a compact anti-self-dual Einstein orbifold with $b_2(\mathcal{M})=k$ and
$\pi^{orb}_1(\mathcal{M})=e$.  Let $n=k+2$.
Then there are distinct real points $t_1,t_2,\ldots,t_{n}\in E$ so that for
$t\in E\setminus\{t_1,t_2,\ldots,t_{n}\}$, $X_t\subset\mathcal{Z}$ is a suborbifold.
Furthermore $X_t$ is a special symmetric toric Fano surface.  The anti-canonical cycle of $X_t$
is $C_1,C_2,\ldots, C_{2n}$, and the corresponding stabilizers are
$\rho_1, \rho_2, \ldots, \rho_{2n}$ which define the vertices in $N=\Z\times\Z$ of
$\Delta^*$ with $X_t =X_{\Delta^*}$.

For $t_i\in E$, $X_{t_i} =D+\ol{D}$, where $D,\ol{D}$ are irreducible degree one divisors with
$\sigma(D)=\ol{D}$.  The $D,\ol{D}$ are suborbifolds of $\mathcal{Z}$ and are toric Fano
surfaces.  We have $D\cap\ol{D}=P_{i}$ and the elements
$\pm(\rho_1,\ldots,\rho_i,-\rho_i +\rho_{i+1},\rho_{n+i+1},\ldots, \rho_{2n})$
define the augmented fans for $D$ and $\ol{D}$.
\end{thm}

Note that a consequence of the theorem is that $\Sigma$ given by equation (\ref{eq:twist-sing}) is the set with
non-trivial stablizers of the action of the complex torus $T^2_\C$.

\section{Sasakian submanifolds}\label{sec:Sasak-submfd}

In this section we use the results on toric 3-Sasakian manifolds to produce a new infinite series of toric
Sasaki-Einstein 5-manifolds corresponding to the toric 3-Sasakian 7-manifolds discussed above and completing
the correspondence in diagram~\ref{int:diag-cor}.
But first we review the Smale/Barden classification of smooth
5-manifolds which will be used.  The possible diffeotypes of toric examples in dimension 5 is
very limited.

\subsection{Classification of 5-Manifolds}
Closed smooth simply connected
spin 5-manifolds were classified by S. Smale~\cite{Sm}.  Subsequently D. Barden extended the classification to the
non-spin case~\cite{Bar}.  Consider the primary decomposition
\begin{equation}\label{H2-decomp}
H_2(M,\Z)\cong \Z^r \oplus \Z_{k_1}\oplus\Z_{k_2}\oplus\cdots\oplus\Z_{k_s},
\end{equation}
where $k_j$ divides $k_{j+1}$.  Of course the decomposition is not unique,
but the $r,k_1,\ldots,k_s$ are.  The second Stiefel-Whitney class defines a homomorphism
$w_2 :H_2(M,\Z)\rightarrow \Z_2$.  One can arrange the decomposition~\ref{H2-decomp} so that $w_2$ is non-zero
on only one component $\Z_{k_j}$, or $\Z$ of $\Z^r$.  Then define $i(M)$ to be $i$ if $2^i$ is the $2$-primary
component of $\Z_{k_j}$, or $\infty$.  Alternatively, $i(M)$ is the minimum $i$ so that $w_2$ is non-zero on a
$2$-primary component of order $2^i$ of $H_2(M,\Z)$.
\begin{thm}[\cite{Sm,Bar}]\label{thm:Sm-Bar}
Smooth simply connected closed 5-manifolds are classifiable up to diffeomorphism.  Any such manifold is
diffeomorphic to one of
\[ X_j \#M_{k_1} \#\cdots\#M_{k_s}, \]
where $-1\leq j\leq\infty, s\geq 0, 1< k_1$ and $K_i$ divides $k_{i+1}$ or $k_{i+1}=\infty$.  A complete set of
invariants is given by $H_2(M,\Z)$ and $i(M)$, and the manifolds $X_{-1},X_0,X_j,X_\infty,M_k,M_\infty$ are
as follows

\begin{center}
\begin{tabular}{|l|c|c|} \hline
$M$ & $H_2(M,\Z)$ & $i(M)$\\ \hline
$M_0=X_0=S^5$ & 0 & 0 \\ \hline
$M_k, 1<k<\infty$ & $\Z_k\oplus\Z_k$ & 0 \\ \hline
$M_\infty =S^2\times S^3$ & $\Z$ & 0 \\ \hline
$X_{-1}$ & $\Z_2$ & 1 \\ \hline
$X_j, 0<j<\infty$ & $\Z_{2^j}\oplus\Z_{2^j}$ & $j$ \\ \hline
$X_\infty$ & $\Z$ & $\infty$ \\ \hline
\end{tabular},\\
\end{center}
where $X_{-1}=SU(3)/SO(3)$ is the Wu manifold, and $X_\infty$ is the non-trivial $S^3$-bundle over
$S^2$.
\end{thm}

The existence of an effective $T^3$ action severely restricts the topology by the following theorem of Oh~\cite{Oh}.
\begin{thm}\label{thm:toric-5mfd}
Let $M$ be a compact simply connected 5-manifold with an effective $T^3$-action.  Then
$H_2(M,\Z)$ has no torsion.  Thus, $M$ is diffeomorphic to
$S^5$, $\#kM_\infty$, or $X_\infty \#(k-1)M_\infty$, where $k=b_2(M)\geq 1$.  Conversely,
these manifolds admit effective $T^3$-actions.
\end{thm}
By direct construction C. Boyer, K. Galicki, and L. Ornea~\cite{BGO} showed that the manifolds in this
theorem admit toric Sasakian structures and, in fact, admit regular Sasakian structures.

A simply connected Sasaki-Einstein manifold must have $w_2=0$, therefore we have the following:
\begin{cor}\label{cor:toric-5mfd-Ein}
Let $M$ be a compact simply connected 5-manifolds with a toric Sasaki-Einstein structure.  Then
$M$ is diffeomorphic to $S^5$, or $\#kM_\infty$, where $k=b_2(M)\geq 1$.
\end{cor}

\subsection{Sasakian Submanifolds and Examples}

Associated to each compact toric anti-self-dual Einstein orbifold $\mathcal{M}$ with
$\pi_1^{orb}(\mathcal{M})=e$ is the twistor space $\mathcal{Z}$ and a family of
embeddings $X_t\subset\mathcal{Z}$ where $t\in E\setminus\{t_1,t_2,\ldots,t_{k+2}\}$
and $X=X_t$ is the symmetric toric Fano surface canonically associated to $\mathcal{M}$.
We denote the family of embeddings by
\begin{equation}
\iota_t: X\rightarrow\mathcal{Z}.
\end{equation}

Let $M$ be the total space of the $S^1$ V-bundle associated to $\mathbf{K}_X$
or $\mathbf{K}_X^{\frac{1}{2}}$, depending on whether $\ind(X)=1$ or $2$.
\begin{thm}\label{thm:subman}
Let $\mathcal{M}$ be a compact toric anti-self-dual Einstein orbifold with
$\pi_1^{orb}(\mathcal{M})=e$.  There exists a Sasakian structure
$\{\tilde{g},\tilde{\Phi},\tilde{\xi},\tilde{\eta}\}$ on
$M$, such that if $(X,\tilde{h})$ is the K\"{a}hler structure making $\pi:M\rightarrow X$ a
Riemannian submersion, then we have the following diagram where the horizontal maps are
isometric embeddings.
\begin{equation}\label{diag:subman-fund}
\beginpicture
\setcoordinatesystem units <1pt, 1pt> point at 0 30
\put {$M$} at -15 60
\put {$\mathcal{S}$} at 15 60
\put {$X$} at -15 30
\put {$\mathcal{Z}$} at 15 30
\put {$\mathcal{M}$} at 15 0
\put {$\ol{\iota}_t$} at 0 65
\put {$\iota_t$} at 0 35
\arrow <2pt> [.3, 1] from -6 60 to 6 60
\arrow <2pt> [.3, 1] from -6 30 to 6 30
\arrow <2pt> [.3, 1] from -15 51 to -15 39
\arrow <2pt> [.3, 1] from 15 51 to 15 39
\arrow <2pt> [.3, 1] from 15 21 to 15 9
\endpicture
\end{equation}

If the 3-Sasakian space $\mathcal{S}$ is smooth, then so is $M$.  If $M$ is smooth, then
\[M\underset{\text{diff}}{\cong}\# kM_\infty, \text{  where  } k=2b_2(\mathcal{S})+1. \]
\end{thm}
\begin{proof}
Let $\{g,\Phi,\xi,\eta\}$ be the fixed Sasakian structure on $\mathcal{S}$ with $\Phi$ descending to the complex structure
on $\mathcal{Z}$.
The adjunction formula gives
\begin{equation}\label{eq:sub-adj}
\mathbf{K}_X \cong\mathbf{K}_\mathcal{Z}\otimes[X]|_X =\mathbf{K}_\mathcal{Z}\otimes\mathbf{K}_\mathcal{Z}^{-\frac{1}{2}}|_X
=\mathbf{K}_\mathcal{Z}^{\frac{1}{2}}|_X.
\end{equation}
Let $h$ be the K\"{a}hler-Einstein metric on $\mathcal{Z}$ related to the 3-Sasakian metric
$g$ on $\mathcal{S}$ by Riemannian submersion.
Recall that $\mathcal{S}$ is the total space of the $S^1$ V-bundle associated to
$\mathbf{L}^{-1}$, (resp. $\mathbf{L}^{-\frac{1}{2}}$ if $w_2(\mathcal{M})=0$).
Also $M$ is the total space of the $S^1$ V-bundle associated to $\mathbf{K}_X$,
(resp. $\mathbf{K}_X^{\frac{1}{2}}$ if $w_2(\mathcal{M})=0$).  Using the isomorphism in
(\ref{eq:sub-adj}) we lift $\iota_t$ to $\ol{\iota}_t$.  The metric on $\mathcal{S}$ is
\[ g=\eta\otimes\eta + \pi^*h.\]
Pull back $h$ and $\eta$ to $\tilde{h}=\iota_t^*h$ and $\tilde{\eta}=\ol{\iota}_t^*\eta$ respectively.
It follows that $\frac{1}{2}\tilde{\eta}=\tilde{\omega}$, where $\tilde{\omega}$ is the K\"{a}hler form of
$\tilde{h}$.  Then $\tilde{\Phi}=\tilde{\nabla}\tilde{\xi}$ is a lift of the complex structure $J$ on $X$.
And the integrability condition in definition~\ref{defn:Sasak} follows from the integrability of $J$.
Then
\[ \tilde{g}=\tilde{\eta}\otimes\tilde{\eta} +\pi^*\tilde{h} \]
is a Sasakian metric on $M$.

If $\mathcal{S}$ is smooth, then the orbifold uniformizing groups act on $\mathbf{L}^{-1}$ (or $\mathbf{L}^{\frac{1}{2}}$)
minus the zero section without non-trivial stabilizers.  By (\ref{eq:sub-adj}) this holds for the bundle
$\mathbf{K}_X$ (or $\mathbf{K}_X^{\frac{1}{2}}$) on $X$.

It follows from $\pi_1^{orb}(\mathcal{M})=e$ that $\pi_1^{orb}(X)=e$ (cf.~\cite{HS2}).
And $\pi_1^{orb}(M)=e$ by arguments as after Proposition~\ref{KE-Sasak}.
If $M$ is smooth, Corollary~\ref{cor:toric-5mfd-Ein} gives the diffeomorphism.
\end{proof}

We are more interested in $M$ with the Sasaki-Einstein metric that exists by corollary
~\ref{cor:Sasak-Einst}.  In this case the horizontal maps are not isometries.

Consider the reducible cases, $t_i\in E,i=1,\ldots, k+2$, where $X_{t_i}=D+\ol{D}\subset\mathcal{Z}$.
Then restricting $\mathbf{L}^{-1}$, (resp. $\mathbf{L}^{-\frac{1}{2}}$ if $w_2(\mathcal{M})=0$), and arguing as above
we obtain Sasakian manifolds $N_i,\ol{N}_i\subset\mathcal{S}$, where smoothness follows from that of $\mathcal{S}$,
whose Sasakian structures $(N_i,g_i,\Phi_i,\xi_i,\eta_i)$ and $(\ol{N}_i,\ol{g}_i,\ol{\Phi}_i,\ol{\xi}_i,\ol{\eta}_i)$
pull back from that of $\mathcal{S}$.
And $N_{t_i}\cap \ol{N}_{t_i}$ is a lens space with the constant curvature metric.
Note that $\varsigma$ restricts to an isometry $\varsigma:N_i \rightarrow\ol{N}_i$ which gives a conjugate isomorphism
between Sasakian structures.  These manifolds do not satisfy the conditions of Proposition~\ref{CY-cond}, so cannot
be transversally deformed to Sasaki-Einstein structures.  But $c_1^B(N_i)>0$, so they can be transversally
deformed to positive Ricci curvature Sasakian by the transverse Calabi-Yau theorem~\cite{EKN}.
By Theorem~\ref{thm:toric-5mfd} $N_i$ is diffeomorphic to $\#kM_\infty$, or $X_\infty \#(k-1)M_\infty$,
where $k=b_2(N_i)=b_2(\mathcal{S})$.  By the remarks after Proposition~\ref{prop:pos-Sasak} $N_i, i=1,\ldots, k+2$,
is diffeomorphic to $X_\infty \#(k-1)M_\infty$ if $w_2(\mathcal{M})\neq 0$.

The family of submanifolds $\ol{\iota}_t :M\rightarrow\mathcal{S}$ for
$t\in E\setminus\{t_1,t_2,\ldots,t_{k+2}\}$ and $N_i,\ol{N}_i,i=1,\ldots,k+2$ for the reducible cases have a simple description.
Recall the 1-form $\eta=\eta^2 +i\eta^3$ of section~\ref{subsec:3Sasak} which is $(1,0)$ with
respect to the CR structure $I=-\Phi^1$.  For $t\in\mathfrak{t}$ let $Y_t$ denote the killing vector field on
$\mathcal{Z}$ with lift $\ol{Y}_t\in I(\mathcal{S},g)$. Then
$\theta(Y_t)\in H^0(\mathcal{Z},\mathcal{O}(\mathbf{L}))$ which defines a holomorphic function on
$\mathbf{L}^{-1}$.  Since $\mathcal{S}$ is the $S^1$ subbundle of $\mathbf{L}^{-1}$, we have $\theta(Y_t)=\eta(\ol{Y}_t)$.
Complexifying gives the same equality for $t\in\mathfrak{t}_{\C}$.
Thus for $t\in E\setminus\{t_1,t_2,\ldots,t_{k+2}\}$, we have
$M_t :=\ol{\iota}_t(M)=(\eta(\ol{Y}_t))\subset\mathcal{S}$ and $N_i\cup\ol{N}_i =(\eta(\ol{Y}_{t_i}))\subset\mathcal{S}$,
where of course $(\eta(\ol{Y}_t))$ denotes the submanifold $\eta(\ol{Y}_t)=0$.
Note that here we are setting 2/3 s of the moment map to zero.

This gives us the new infinite families of Sasaki-Einstein manifolds and the diagram~\ref{int:diag-cor}.
\begin{thm}\label{thm:main}
Let $(\mathcal{S},g)$ be a toric 3-Sasakian 7-manifold with $\pi_1(\mathcal{S})=e$.
Canonically associated to $(\mathcal{S},g)$ are a special symmetric toric Fano surface $X$ and a
toric  quasi-regular Sasaki-Einstein 5-manifold $(M,g,\Phi,\xi,\eta)$ which fit in the commutative diagram
~\ref{diag:subman-fund}.
We have $\pi_1^{orb}(X)=e$ and $\pi_1(M)=e$.  And
\[M\underset{\text{diff}}{\cong}\# kM_\infty, \text{  where  } k=2b_2(\mathcal{S})+1 \]
Furthermore $(\mathcal{S},g)$ can be recovered from either $X$ or $M$ with their torus actions.
\end{thm}
\begin{proof}
 Note that the homotopy sequence
 \[\cdots\rightarrow\pi_1(G)\rightarrow\pi_1(\mathcal{S})\rightarrow\pi^{orb}_1(\mathcal{M})\rightarrow e,\]
where $G=SO(3)$ or $Sp(1)$, shows that $\pi_1^{orb}(\mathcal{M})=e$.  The toric surface $X_{\Delta^*}$ and
Sasakian 5-manifold $(M,\tilde{g},\tilde{\Phi},\tilde{\xi},\tilde{\eta})$ are given in Theorem~\ref{thm:subman}.
By Corollary~\ref{cor:Sasak-Einst} this Sasakian structure has a transversal deformation to a Sasaki-Einstein
structure.

Given either $X$ or $M$, by the discussion in section~\ref{sec:asd-Einst}, the orbifold $\mathcal{M}$
can be recovered.  By Theorem~\ref{thm:asd-class} $\mathcal{M}$ admits
a unique anti-self-dual Einstein structure up to homothety compatible with the torus action.  
The 3-Sasakian space $\mathcal{S}$ and its twistor space $\mathcal{Z}$ can be constructed from $\mathcal{M}$ 
and its anti-self-dual Einstein structure. (cf.~\cite{BG2}).
\end{proof}

\begin{cor}\label{cor:inf-xpl}
For each odd $k\geq 3$ there is a countably infinite number of distinct toric quasi-regular Sasaki-Einstein
structures on $\# kM_\infty$.
\end{cor}

We do not know the Sasaki-Einstein metrics explicitly.  But if $c=\ind(M)$, an application of
Corollary~\ref{cor:toric-fano-vol} gives
\[\begin{split}
   \Vol(M,g) &=\frac{2\pi c}{3}\Vol(X,\omega)\\
         &=2c\left(\frac{\pi}{3}\right)^3\Vol(\Sigma_{-k}),\\
  \end{split}\]
where $\omega$ is the transversal K\"{a}hler metric.
We have $c=1$ or $2$.  Let $(M,g_i), i\in\Z_+,$ be any infinite sequence of metrics on $\# kM_\infty$ in Corollary
~\ref{cor:inf-xpl}.
These Sasaki-Einstein structures have leaf spaces
$X_i$, where $X_i =X_{\Delta^*_i}$.  Observe that the polygons $\Delta^*_i$ get arbitrarily large,
and the anti-canonical polytopes $(\Sigma_{-k})_i$ satisfy
\[ \Vol((\Sigma_{-k})_i)\rightarrow 0,\text{  as  }i\rightarrow\infty.\]
Thus we have $\Vol(M,g_i)\rightarrow 0$, as $i\rightarrow\infty$.

\subsection{Examples}

We consider some of the examples obtained starting with the simplest.  In particular
we can determine some of the spaces in diagram \ref{int:diag-cor} associated to a
toric 3-Sasakian 7-manifold more explicitly in some cases.

\subsubsection{Smooth examples}

It is well known that there exists only two complete examples of positive scalar curvature
anti-self-dual Einstein manifolds~\cite{Hi2}~\cite{FK}, $S^4$ and $\cps^2$ with the round and Fubini-Study metrics respectively.
Note that we are considering $\cps^2$ with the opposite of the usual orientation.
\\

\noindent
\emph{$\mathcal{M}= S^4$}

For the spaces in diagram~\ref{int:diag-cor} we have:
$\mathcal{M}= S^4$ with the round metric; its twistor space $\mathcal{Z}=\cps^3$ with the Fubini-study metric;
the quadratic divisor $X\subset\mathcal{Z}$ is $\cps^1\times\cps^1$ with the
homogeneous K\"{a}hler-Einstein metric; $M = S^2\times S^3$ with the homogeneous Sasaki-Einstein structure;
and $\mathcal{S}= S^7$ has the round metric.
In this case diagram~\ref{int:diag-cor} becomes the following.

\begin{equation}
\beginpicture
\setcoordinatesystem units <1pt, 1pt> point at 0 30
\put {$S^2\times S^3$} at -30 60
\put {$S^7$} at 30 60
\put {$\cps^1\times\cps^1$} at -30 30
\put {$\cps^3$} at 30 30
\put {$S^4$} at 30 0
\arrow <2pt> [.3, 1] from -1 60 to 11 60
\arrow <2pt> [.3, 1] from -1 30 to 11 30
\arrow <2pt> [.3, 1] from -30 51 to -30 39
\arrow <2pt> [.3, 1] from 30 51 to 30 39
\arrow <2pt> [.3, 1] from 30 21 to 30 9
\endpicture
\end{equation}

This is the only example, I am aware of, for which the horizontal maps are isometric immersions
when the toric surface and Sasakian space are equipped with the Einstein metrics.
\\

\noindent
\emph{$\mathcal{M}=\cps^2$}

In this case $\mathcal{M}=\cps^2$ with the Fubini-Study metric and the reverse of the usual orientation; its twistor
space is $\mathcal{Z}=F_{1,2}$, the manifold of flags $V\subset W\subset\C^3$ with $\dim V =1$ and $\dim W=2$,
with the homogeneous K\"{a}hler-Einstein metric.
The projection $\pi:F_{1,2}\rightarrow\cps^2$ is as follows.  If $(p,l)\in F_{1,2}$ so
$l$ is a line in $\cps^2$ and $p\in l$, then $\pi(p,l)=p^\perp \cap l$, where $p^\perp$ is the
orthogonal compliment with respect to the standard hermitian inner product.
We can define $F_{1,2}\subset\cps^2\times(\cps^2)^*$ by
\[ F_{1,2}=\{([p_0:p_1:p_2],[q^0:q^1:q^2])\in\cps^2\times(\cps^2)^*: \sum p_iq^i =0\}.\]
And the complex contact structure is given by $\theta=q^idp_i -p_idq^i$.
Fix the action of $T^2$ on $\cps^2$ by
\[ (e^{i\theta},e^{i\phi})[z_0:z_1:z_2]=[z_0:e^{i\theta}z_1:e^{i\phi}z_2].\]
Then this induces the action on $F_{1,2}$
\[ (e^{i\theta},e^{i\phi})([p_0:p_1:p_2],[q^0:q^1:q^2])=([p_0:e^{i\theta}p_1:e^{i\phi}p_2],[q^0:e^{-i\theta}q^1:e^{-i\phi}q^2]). \]
Given $[a,b]\in\cps^1$ the one parameter group $(e^{ia\tau},e^{ib\tau})$ induces the holomorphic vector field
$W_\tau\in\Gamma(T^{1,0}F_{1,2})$ and the quadratic divisor $X_\tau=(\theta(W_\tau))$ given by
\[ X_\tau =(ap_1q^1+bp_2q^2=0,\quad p_iq^i=0).\]
One can check directly that $X_\tau$ is smooth for $\tau\in\cps^1\setminus\{[1,0],[0,1],[1,1]\}$ and
$X_\tau =\cps^2_{(3)}$, the equivariant blow-up of $\cps^2$ at 3 points.  For
$\tau\in\{[1,0],[0,1],[1,1]\}$, $X_\tau =D_\tau +\ol{D}_\tau$ where
both $D_\tau,\ol{D}_\tau$ are isomorphic to the Hirzebruch surface
$F_1 =\mathbb{P}(\mathcal{O}_{\cps^1}\oplus\mathcal{O}_{\cps^1}(1))$.

The Sasaki-Einstein space is $M=\#3(S^2\times S^3)$.
And the Sasakian manifolds $N_\tau,\ol{N}_\tau$ are diffeomorphic to $X_\infty$.
$\mathcal{S}=\mathcal{S}(1,1,1)=SU(3)/U(1)$ with the homogeneous 3-Sasakian structure.
This case has the following diagram.

\begin{equation}
\beginpicture
\setcoordinatesystem units <1pt, 1pt> point at 0 30
\put {$\#3(S^2\times S^3)$} at -35 60
\put {$SU(3)/U(1)$} at 35 60
\put {$\cps^2_{(3)}$} at -35 30
\put {$F_{1,2}$} at 35 30
\put {$\cps^2$} at 35 0
\arrow <2pt> [.3, 1] from -6 60 to 6 60
\arrow <2pt> [.3, 1] from -6 30 to 6 30
\arrow <2pt> [.3, 1] from -35 51 to -35 39
\arrow <2pt> [.3, 1] from 35 51 to 35 39
\arrow <2pt> [.3, 1] from 35 21 to 35 9
\endpicture
\end{equation}

\subsubsection{Galicki-Lawson quotients}

The simplest examples of quaternionic-K\"{a}hler quotients are the Galicki-Lawson examples first appearing in~\cite{GL}
and further considered in~\cite{BGMR}.
These are circle quotients of $\qps^2$.  In this case the weight matrices are of the form
$\Omega=\mathbf{p}=(p_1,p_2,p_3)$ with the admissible set
\[\{\mathcal{A}_{1,3}(Z)=\{\mathbf{p}\in\Z^3| p_i\neq 0 \text{ for }i=1,2,3 \text{ and } \gcd(p_i,p_j)=1\text{ for }i\neq j\} \]
We may take $p_i>0$ for $i=1,2,3$.
The zero locus of the 3-Sasakian moment map $N(\mathbf{p})\subset S^{11}$ is diffeomorphic to the
Stiefel manifold $V_{2,3}^\C$ of complex 2-frames in $\C^3$ which can be identified as
$V_{2,3}^\C\cong U(3)/U(1)\cong SU(3)$.
Let $f_{\mathbf{p}}: U(1)\rightarrow U(3)$ be
\[ f_{\mathbf{p}}(\tau)=\begin{bmatrix}\tau^{p_1} & 0 & 0 \\ 0 & \tau^{p_2} & 0 \\ 0 & 0 & \tau^{p_3}\end{bmatrix}. \]
Then the 3-Sasakian space $\mathcal{S}(\mathbf{p})$ is diffeomorphic to the quotient of $SU(3)$ by the
action of $U(1)$
\[\tau\cdot W=f_{\mathbf{p}}(\tau)Wf_{(0,0,-p_1 -p_2 -p_3)}(\tau)\text{  where  }\tau\in U(1)\text{ and }W\in SU(3). \]
Thus $\mathcal{S}(\mathbf{p})\cong SU(3)/U(1)$ is a biquotient similar to the examples considered by
Eschenburg in~\cite{Es}.

The action of the group $SU(2)$ generated by $\{\xi^1,\xi^2,\xi^3\}$ on
$N(\mathbf{p})\cong SU(3)$ commutes with the action of $U(1)$.
We have $N(\mathbf{p})/SU(2)\cong SU(3)/SU(2)\cong S^5$ with $U(1)$ acting by
\[\tau\cdot v=f_{(-p_2-p_3,-p_1-p_3,-p_1-p_2)}v \text{  for  }v\in S^5\subset\C^3. \]
We see that $\mathcal{M}_{\Omega}\cong\cps^2_{a_1,a_2,a_3}$ where $a_1=p_2+p_3, a_2=p_1+p_3, a_3=p_1+p_2$
and the quotient metric is anti-self-dual with the reverse of the usual orientation.
If $p_1,p_2,p_3$ are all odd then the generic leaf of the 3-Sasakian foliation $\mathscr{F}_3$ is
$SO(3)$.  If exactly one is even, then the generic leaf is $Sp(1)$.
Denote by $X_{p_1,p_2,p_3}$ the toric Fano divisor, which can be considered as a generalization of
$\cps^2_{(3)}$.
We have the following spaces and embeddings.

\begin{equation}
\beginpicture
\setcoordinatesystem units <1pt, 1pt> point at 0 30
\put {$\#3(S^2\times S^3)$} at -38 60
\put {$\mathcal{S}(p_1,p_2,p_3)$} at 35 60
\put {$X_{p_1,p_2,p_3}$} at -30 30
\put {$\mathcal{Z}(p_1,p_2,p_3)$} at 35 30
\put {$\cps^2_{a_1,a_2,a_3}$} at 45 0
\arrow <2pt> [.3, 1] from -8 60 to 4 60
\arrow <2pt> [.3, 1] from -10 30 to 2 30
\arrow <2pt> [.3, 1] from -30 51 to -30 39
\arrow <2pt> [.3, 1] from 35 51 to 35 39
\arrow <2pt> [.3, 1] from 35 21 to 35 9
\endpicture
\end{equation}

\section{Positive Ricci curvature examples}\label{sec:pos-Sasak}

In this section we use the toric geometry developed to construct examples of positive Ricci curvature
Sasakian structures on the manifolds $X_\infty \#(k-1)M_\infty$ in Theorem~\ref{thm:Sm-Bar}.  By Theorem
~\ref{thm:toric-5mfd} these are the only simply connected non-spin 5-manifolds that can admit toric Sasakian
structures.  They are already known to admit Sasakian structures~\cite{BGO}.

Define a marked fan $\Delta^* =\Delta^*_{k,p}, k\geq 2, p\geq 0$ as follows.
Let $\sigma_0 = (-1,0)$, $\sigma_1=(0,1);$
$\sigma_j=(j-1, \frac{j(j-1)}{2}-1),j=2,\ldots,k;$
$\sigma_{k+1}=(k,\frac{(k+1)k}{2}-1+p)$, $\sigma_{k+2}=(0,\frac{(k+1)k}{2}+p)$.
For $k=1$ define $\Delta^*_{1,p},p\geq 0$ by $\sigma_0=(-1,0),\sigma_1 =(0,1),
\sigma_2=(1,1+p),\sigma_3=(0,2+p)$.
And define $l\in\SF(\Delta^*)$ by
\[ l(\sigma_j):=\begin{cases}
 0 & \text{if $j=0$,}\\
 -1 & \text{if $1\leq j\leq k+2$}.
\end{cases}\]
It is easy to check that $l$ is strictly upper convex, and we have $\pi_1^{orb}(X)=e$ for all of the above fans.
The cone $\mathcal{C}_l$ corresponding to $(\mathbf{L}_l^{-1})^\times$ on $X_{\Delta^*}$ satisfies
the smoothness condition of equation~\ref{eq:toric-smooth}, and $(\mathbf{L}_l^{-1})^\times$ is simply connected.
Also, $-k\in\SF(\Delta^*)$ is strictly upper convex.  Thus by Corollary~\ref{cor:toric-fano} $2\pi c_1(X_{\Delta^*})=[\omega]$ for
a K\"{a}hler form $\omega$.
\begin{prop}\label{prop:pos-Sasak}
Each of the manifolds $M_{k,p}$, for $k\geq 1,p\geq 0$, is diffeomorphic to $X_\infty \#(k-1)M_\infty$.
And for each $p$ has a distinct positive Ricci curvature Sasakian structure.
\end{prop}
\begin{proof}
Let $X =X_{\Delta^*_{k,p}}$, and let $\mathbf{L}$ be the holomorphic line V-bundle associated with
$l\in\SF(\Delta^*)$.  Since $l$ is strictly upper convex, if $\omega$ is the canonical metric on $X$ of
$\Sigma^*_l$, then $[\omega]\in 2\pi c_1(\mathbf{L})$ by Theorem~\ref{thm:toric-ampl}.
By Proposition~\ref{KE-Sasak} there is a Sasakian structure $(M,g,\Phi,\xi, \eta)$
where $M$ is the principle $S^1$ subbundle of $\mathbf{L}^{-1}$.
The the toric cone $\mathcal{C}$ of $(\mathbf{L}^{-1})^\times \cong C(M)$ satisfies the smoothness condition.
Also, by Theorem~\ref{thm:toric-ampl} $\mathbf{K}^{-1} >0$.  In other words $c_1^B(M)$ has a positive representative.
By the transverse Calabi-Yau theorem~\cite{EKN} there is a transversal deformation $(M,\tilde{g},\tilde{\Phi},\xi,\tilde{\eta})$
with $\tilde{\eta}=\eta +2d^c \phi$, for $\phi\in C^\infty_B(M)$, with $\Ric^T_{\tilde{g}}$ positive.
Then by (\ref{eq:Sasak-Ric}) a homothetic deformation $\tilde{g}_a =a^2\tilde{\eta}\otimes\tilde{\eta} + a\tilde{g}^T$ for small enough
$a\in\R_+$ has $\Ric_{\tilde{g}_a}$ positive on $M$.

Considering $\pi:M\rightarrow X$ as an $S^1$ V-bundle over $X$, it lifts to a genuine fiber bundle over $B(X)=M\times_{S^1} E(S^1)$,
$\tilde{\pi}:\tilde{M}\rightarrow B(X)$.  Here $E(S^1)$ is the universal $S^1$-principal bundle, and $B(X)$
is the orbifold classifying space (cf.~\cite{BG3}).  Since $M$ is smooth, we have a homotopy equivalence
$\tilde{M}\simeq M$.  And since $\pi_1^{org}(X)=e$, the first few terms of the Gysin sequence give
\[ \Z\overset{\cup e}{\longrightarrow} H^2_{orb}(X,\Z)\overset{\pi^*}{\longrightarrow} H^2(M,\Z)\rightarrow 0, \]
where $e=c_1^{orb}(\mathbf{L})\in H^2_{orb}(X,\Z)$.
We have $w_2(M)\equiv c_1(D) \mod 2\equiv \pi^*(c_1^{orb}(X)) \mod 2$,
which is zero precisely when $x=c_1^{orb}(X)$ is divisible by 2 in $H^2_{orb}(X,\Z)/\Z(e)$.
If this is the case then there is a $u\in H^2_{orb}(X,\Z)$ with
$2u = x + al$, with $a\in\Z$ odd.  It is easy to see that $u\in H^2_{orb}(X,\Z)$ is represented by a holomorphic
line V-bundle, since $2u$ is.  By Proposition~\ref{prop:toric-vbund}, this V-bundle is $\mathbf{L}_u$ for some
$u\in\SF(\Delta^*)$, and we have an equation
\begin{equation}\label{eq:sf}
 2u =-k +al +f, \quad \text{ where } f\in M.
\end{equation}
Evaluating (\ref{eq:sf}) on $\sigma_0$ gives $2u(\sigma_0)=-1 +f(\sigma_0)$, and on $\sigma_2 =-\sigma_0$ gives
$2u(\sigma_2) =-(1+a) -f(\sigma_0)$.  The first equation implies $f(\sigma_0)$ is odd, the second that $f(\sigma_0)$
is even.  Thus $w_2(M)\neq 0$. Theorem~\ref{thm:toric-5mfd} then completes the proof.
\end{proof}

One can use arguments of M. Demazure (cf.~\cite{O2} \S 3.4) as in the smooth case to show that $\Aut^o(X)$ is not reductive
for these examples.  The fan $\Delta_{k,p}^*$ has two roots $\alpha_1,\alpha_2$ and $\alpha_2\neq -\alpha_1$
as the Hirzebruch surface $F_1$.  The Lie algebra of Hamiltonian holomorphic vector fields of $M$ is not reductive.
By the proof of the Lichnerowicz theorem in the Sasakian case~\cite{BGS} the Sasakian structure can not be transversally deformed to
constant scalar curvature.

Let $N_i,\ol{N}_i\subset\mathcal{S}, i=1,\ldots, k+2$ be the Sasakian submanifolds of the 3-Sasakian manifold
$\mathcal{S}$ as discussed after Theorem~\ref{thm:subman}.  Suppose $w_2(\mathcal{M})\neq 0$.  So $N_i,\ol{N}_i$
are principle $S^1$ subbundles of $\mathbf{L}^{-1}$ restricted to toric surfaces $D,\ol{D}$ with fans as in
Theorem~\ref{thm:twistor-divisor}.  The augmented fan $\Delta^*$ of $D$ has elements
$\sigma_1 =\rho_1,\ldots,\sigma_i =\rho_i,\sigma_{i+1}=-\rho_i +\rho_{i+1},\sigma_{i+2}=\rho_{n+i+1},\ldots,
\sigma_{k+3}=\rho_{2n}$, $n=k+2$, and $\mathbf{L}|_D$ is the line bundle associated to $l\in\SF(\Delta^*)$
with
\[ l(\sigma_j):=\begin{cases}
 0 & \text{if $j=i+1$,}\\
 -1 & \text{otherwise}.
\end{cases}\]
The same argument in the above proposition shows that $N_i$ is diffeomorphic to $X_\infty \#(k-1)M_\infty$.

\section{Higher dimensional examples}\label{sec:high-dim}

In this section we employ the join construction of C. Boyer and K. Galicki~\cite{BGO,BG3} to construct higher
dimensional examples.  Let $(M_i,g_i,\Phi_i,\xi_i,\eta_i), i=1,2$ be quasi-regular Sasakian manifolds of dimensions
$2m_i +1, i=1,2$..
Make homothetic deformations of the Sasakian structures so that the $S^1$-actions generated by $\xi_i, i=1,2$ have
period $1$.  Then the transverse K\"{a}hler forms $\omega_i^T,i=1,2$ descend to forms $\omega_i,i=1,2$ on the leaf spaces
$\mathcal{Z}_i,i=1,2$ with $[\omega_i]\in H^2_{orb}(\mathcal{Z}_i,\Z),i=1,2$.  Then for a pair of positive integers
$(k_1,k_2)$ we have a k\"{a}hler form $k_1\omega_1 +k_2\omega_2$ on $\mathcal{Z}_1\times\mathcal{Z}_2$ with
$[k_1\omega_1 +k_2\omega_2]\in H^2_{orb}(\mathcal{Z}_1\times\mathcal{Z}_2,\Z)$.
By Proposition~\ref{KE-Sasak} there is an $S^1$ V-bundle, denoted $M_1\star_{k_1,k_2}M_2$, with a homothetic family
of Sasakian structures.  We may assume that $\gcd(k_1,k_2)=1$; for if $(k_1,k_2)=(lk_1',lk_2')$,
$M_1\star_{k_1,k_2}M_2 =(M_1\star_{k_1',k_2'}M_2)/\Z_l$.
Note that
\[M_1\star_{k_1,k_2}M_2 = (M_1 \times M_2)/{S^1(k_1,k_2)},\]
where $S^1$ acts by $(x,y)\rightarrow (e^{ik_2\theta}x,e^{-ik_1\theta})$.

In general the join $M_1\star_{k_1,k_2}M_2$ is an orbifold of dimension $2(m_1+m_2)+1$.
But a simple condition exists that implies smoothness.
Let $v_i=\ord(M_i),i=1,2$ denote the $\lcm$ of the orders of the leaf holonomy groups of the $\mathscr{F}_{\xi_i}$.
\begin{prop}[\cite{BG3,BGO}]\label{prop:join-sm}
 For each pair $(k_1,k_2)$ of relatively prime positive integers, $M_1\star_{k_1,k_2}M_2$ is a smooth quasi-regular
 Sasakian manifold if, and only if, $\gcd(v_1 k_2, v_2 k_1)=1$.
\end{prop}

Note that if $M_1$ and $M_2$ are positive, that is, $\ric^T_{g_i}, i=1,2$ are positive, then
$M_1\star_{k_1,k_2}M_2$ is positive.  Suppose $M_1$ and $M_2$ are Sasakian Einstein.
We have $\ind(\mathcal{Z}_1\times\mathcal{Z}_2)=\gcd(\ind(\mathcal{Z}_1),\ind(\mathcal{Z}_2))$.
So we define the \emph{relative indices} of $M_1$ and $M_2$ to be
\begin{equation}
l_i =\frac{\ind(M_i)}{\gcd(\ind(M_1),\ind(M_2))},\quad \text{ for }i=1,2.
\end{equation}
Then the homothetic family of Sasakian structures on the join $M_1\star_{l_1,l_2}M_2$ has a
Sasakian Einstein structure with transverse metric
\[ g^T=\frac{(m_1 +1)g_1^T +(m_2 +1)g_2^T}{m_1 +m_2 +1}.\]

Let $M_k \cong \#k(S^2\times S^3)$ be one of the Sasaki-Einstein manifolds constructed in
section~\ref{sec:Sasak-submfd}.  Then $\ind(M_k)=1$ or $2$.  Consider $S^{2m+1}$ with its standard Sasakian
structure.  Then $\ind(S^{2m+1})=m+1$.  Then the relative indices of $S^{4j+3}$ and $M_k$ are both
$1$.  So Proposition~\ref{prop:join-sm} implies that $M_k\star S^{4j+3}$, $j\geq 0$, is a Sasaki-Einstein
$(4j+7)$-manifold for all examples $M_k$.   One can iterate this procedure; for example
\[ M_k\star\overset{\text{p times}}{\overbrace{S^3\star \cdots\star S^3}},\]
is a $5+2p$ dimensional Sasaki-Einstein manifold.
By making repeated joins to the examples in Corollary~\ref{cor:inf-xpl} we obtain the following.

\begin{prop}
 For every possible dimension $n=2m+1\geq 5$ there are infinitely many toric quasi-regular Sasaki-Einstein
 manifolds with arbitrarily high second Betti number.

 More precisely, in dimension $5$ we have examples of every
 odd $b_2$, and infinitely many examples for each odd $b_2 \geq 3$.

 In dimension $n=2m+1\geq 7$ for $m$ odd we have examples of every even $b_2$, and infinitely many examples for
 each even $b_2 \geq 4$.

 In dimension $n=2m+1\geq 9$ for $m$ even we have examples of every odd $b_2$, and infinitely many examples for
 each odd $b_2 \geq 5$.
\end{prop}

\begin{xpl}
Consider the $7$-dimensional Sasaki-Einstein manifolds $S^3\star M_k$. The Gysin
sequence of $\pi: S^3\star M_k\rightarrow\cps^1 \times X$ determines the cohomology in
$\Q$-coefficients.  The Leray spectral sequence of the fiber bundle $\varpi :S^3\star M_k\rightarrow\cps^1$
with fiber $M_k$ can be used to show the following.
\\
\renewcommand{\arraystretch}{1.5}
\begin{center}
\begin{tabular}{|c|cccccccc|}\hline
$p$ & $0$ & $1$ & $2$ & $3$ & $4$ & $5$ & $6$ & $7$ \\ \hline
$H^p(S^3\star M_k,\Z)$ & $\Z$ & $0$ & $\Z^{k+1}$ & $0$ & $T$ & $\Z^{k+1}$ & $0$ & $\Z$ \\ \hline
\end{tabular}
\end{center}
\vspace{12pt}
Here $T$ is a torsion group.  This is quite similar to the cohomology of $\mathcal{S}_\Omega$
in Theorem~\ref{thm:3-Sasak-coh}.
\end{xpl}

Let $N=N_{k,p}$ denote the positive Ricci curvature Sasakian 5-manifold of Proposition~\ref{prop:pos-Sasak}
diffeomorphic to $X_\infty \#(k-1)M_\infty$.  Then we see by Proposition~\ref{prop:join-sm} that
for any regular Sasaki-Einstein manifold $M$, $N\star_{1,l} M$ is smooth for any $l\geq 1$.
Furthermore, we may transversally deform it to a positive Ricci curvature structure.

Let $\pi_i:\mathbf{L}_i\rightarrow Z_i, i=1,2$ be the holomorphic V-bundles whose $S^1$-principal bundles are
$N$ and $M$ respectively.  Then $N\star_{1,l} M$ is the principal $S^1$ V-bundle of
$\pi :\mathbf{L}_1\otimes\mathbf{L}_2^{l}$, and as in the proof of Proposition~\ref{prop:pos-Sasak} we have
\[ \Z\overset{\cup e}{\longrightarrow} H^2_{orb}(Z_1\times Z_2,\Z)\overset{\pi^*}{\longrightarrow}
H^2(N\star_{1,l} M,\Z)\rightarrow 0, \]
where $e=c_1(\mathbf{L}_1)+lc_1(\mathbf{L}_2)$.  We have
$w_2(N\star_{1,l} M)\equiv\pi^*(c_1(Z_1)+c_1(Z_2)) \mod 2$, and $w_2(N\star_{1,l} M)=0$ if,
and only if, $\pi^*(c_1(Z_1)+c_1(Z_2))$ is divisible by 2 in $H^2(N\star_{1,l} M,\Z)$.  If this is the case then
there is an $u\in H^2_{orb}(Z_1\times Z_2,\Z)$ with
\[ 2u=c_1(Z_1)+c_1(Z_2) +s(c_1(\mathbf{L}_1)+lc_1(\mathbf{L}_2)),\quad\text{for }s\in\Z.\]
Let $\iota: Z_1\rightarrow Z_1\times Z_2$ be the inclusion $\iota(x)=(x,y)$ with $y\in Z_2$ a smooth point.
Then $2\iota^*u =c_1(Z_1)+sc_1(\mathbf{L}_1)$ which contradicts Proposition~\ref{prop:pos-Sasak}.
Thus $w_2(N\star_{1,l} M)\neq 0$.  By taking joins of the $N_{k,p},k\geq 1,p\geq 0$ with spheres as above
we obtain the following.
\begin{prop}
For every dimension $n=2m+1\geq 5$ there exist infinitely many toric quasi-regular positive Ricci
curvature Sasakian manifolds of each $b_2 \geq 2$, and of each $b_2\geq 1$ in dimension $5$.
These examples are simply connected and have $w_2\neq 0$.  Therefore, they do not admit a
Sasaki-Einstein structure.
\end{prop}

\bibliographystyle{plain}

\end{document}